\documentclass[a4paper,12pt]{amsart}
\usepackage{amssymb}

\newtheorem{thm}{Theorem}[section]
\newtheorem{lem}{Lemma}[section]
\newtheorem{prop}{Proposition}[section]
\newtheorem{cor}{Corollary}[section]

\newcommand{\Z}{\mathbf{Z}}
\newcommand{\Q}{\mathbf{Q}}

\def\c{\circeq}
\def\d{\doteq}

\begin{document}

\title[Second homologies of the automorphism group of a free group]
{Twisted second homology groups of the automorphism group of a free group}
\author[Takao Satoh]{Takao Satoh}
\address{Department of Mathematices, Faculty of Science Division II, Tokyo University of Science, \\
     1-3, Kagurazaka, Shinjuku-ku, Tokyo 162-8601, Japan}
\email{takao@rs.tus.ac.jp}

\maketitle
\begin{center}
{\footnotesize Department of Mathematics, Faculty of Science Division II, Tokyo University of Science, \\
         1-3, Kagurazaka, Shinjuku-ku, Tokyo 162-8601, Japan}
\end{center}

\begin{abstract}
In this paper, we compute the second homology groups of the automorphism group of a free group
with coefficients in the abelianization of the free group and its dual group except for the $2$-torsion part, using
combinatorial group theory.
\end{abstract}

\section{Introduction}\label{1}

Let $F_n$ be a free group of rank $n$, and let $\mathrm{Aut}\,F_n$ denote the automorphism group of $F_n$.
There are several remarkable computation of the (co)homology groups of ${\mathrm{Aut}\, F_n}$
with trivial coefficients.
For example, Gersten \cite{Ger} showed $H_2(\mathrm{Aut}$ $\,F_n,\Z)=\Z/2\Z$ for $n \geq 5$, and
Hatcher and Vogtmann \cite{Hat} showed $H_q(\mathrm{Aut}\,F_n,\Q)=0$ for $n \geq 1$ and $1 \leq q \leq 6$, except for
$H_4(\mathrm{Aut} \,F_4, \Q)= \Q$.
In this paper we consider twisted (co)homology groups of ${\mathrm{Aut}\, F_n}$.
Let $H$ be the abelianization of $F_n$ and $H^*:=\mathrm{Hom}_{\Z}(H,\Z)$ the dual group of $H$.
The group ${\mathrm{Aut}\, F_n}$ naturally acts on $H$ and $H^*$.
The main interest of this paper is to compute the homology groups of ${\mathrm{Aut}\, F_n}$
with coefficients in $H$ and $H^*$ using combinatorial groups theory, in particular using
a finite presentation for ${\mathrm{Aut}\, F_n}$.

In our previous paper \cite{Sat}, we computed the first homology groups
of ${\mathrm{Aut}\, F_n}$ with coefficients in $H$ and $H^*$ for $n \geq 2$.
In this paper, we show that the second homology groups $H_2(\mathrm{Aut}\,F_{n},H)$ and $H_2(\mathrm{Aut}\,F_{n},H^*)$
are trivial except for the $\Z/2\Z$-part for $n \geq 6$.
Let $L$ be the subring $\Z[\frac{1}{2}]$ of $\Q$ which is obtained from the ring $\Z$ by attaching $1/2$.
The ring $L$ is a principal ideal domain in which the element $2$ is invertible.
For any $\Z$-module $M$, we denote by $M_L$ the $L$-module $M {\otimes}_{\Z} L$. Then our main theorem is
\begin{thm}\label{X1}
For $n \geq 6$,
\[ H_2(\mathrm{Aut}\,F_n,H_L)=0, \hspace{2em} H_2(\mathrm{Aut}\,F_n,H_L^*)=0. \]
\end{thm}
Recently, Hatcher and Wahl \cite{Nat} showed $H_i(\mathrm{Aut}\,F_{n},H)=0$ for $n \geq 3i+9$ using
the stability of the homology groups of the mapping class groups of certain 3-manifolds.
If $n \geq 15$, one of our result $H_2(\mathrm{Aut}\,F_n,H_L)=0$ is immediately follows from
the results of them.
Our computation, however, is based on combinatorial group theory, and quite different from that of them.
Furthermore we remark that the computation of $H_2(\mathrm{Aut}\,F_n,H_L^*)=0$, to which we cannot apply their method directly,
is more complicated than that of $H_2(\mathrm{Aut}\,F_n,H_L)=0$.

Here we summarize the proof of Theorem {\rmfamily \ref{X1}}. To begin with,
we review the computation of $H_1(\mathrm{Aut}\,F_{n},H)=0$ and $H_1(\mathrm{Aut}\,F_{n},H^*)=\Z$ for $n \geq 4$ due to
\cite{Sat}.
Let $\mathrm{Aut}^{+}F_n$ be the index-$2$ subgroup of $\mathrm{Aut}\,F_n$
defined by the kernel of the composition map of
a natural map $\rho : \mathrm{Aut}\,F_n \twoheadrightarrow \mathrm{Aut}\,H = GL(n,\Z)$
and the determinant map $\mathrm{det} : GL(n,\Z) \longrightarrow \{\pm 1\}$.
To compute the first homology groups of $\mathrm{Aut}\,F_{n}$,
we computed those of $\mathrm{Aut}^{+}F_n$ using a finite presentation
for it due to Gersten \cite{Ger}.
Observing the Lyndon-Hochshild-Serre spectral sequence of
\begin{equation}
 1 \rightarrow \mathrm{Aut^{+}}F_n \rightarrow \mathrm{Aut}\,F_n \rightarrow \{\pm1\} \rightarrow 1, \label{eq200}
\end{equation}
we obtain the results for $\mathrm{Aut}\,F_n$.
We also computed the homology groups for $n=2$ and $3$, and showed they have non-trivial $2$-torsions. For details, see \cite{Sat}.

Now, we show the outline of the computation of the second homology groups.
First, we compute the second homology groups of
$\mathrm{Aut}^{+}F_n$, using a reduced finite presentation $\langle X \,|\, R \rangle$ for it
introduced in Section {\rmfamily \ref{2}}, which is obtained from the Gersten's presentation using Tietze transformations.
Let $F$ and $\bar{R}$ be the free group on $X$ and the normal closure of $R$ in $F$
respectively.
Then, for $M=H_L$ and $H^*_L$, we have a five-term exact sequence
\[\begin{split}
   H_2(F,M) \rightarrow H_2(\mathrm{Aut}^{+}F_n, & M) \rightarrow H_1(\bar{R},M)_{\mathrm{Aut}^{+}F_n} \\
                     & \xrightarrow{\varphi} H_1(F,M) \rightarrow H_1(\mathrm{Aut}^{+}F_n,M) \rightarrow 0
  \end{split}\]
of $L$-modules. Since $F$ is a free group, $H_2(F,M)=0$. Furthermore,
we see $H_1(F,M) = L^{\oplus \{2n(n^2-n)-n \}}$, and we have obtained the rank $r_M$ of
the free $L$-module $H_1(\mathrm{Aut}^{+}F_n,M)$ by the results of \cite{Sat}.
In Section {\rmfamily \ref{3}}, we show that
\[ H_1(\bar{R},M)_{\mathrm{Aut}^{+}F_n} = L^{\oplus \{2n(n^2-n)-n - r_M \}} \]
by reducing generators of $H_1(\bar{R},M)_{\mathrm{Aut}^{+}F_n}$.
This implies that the map $\varphi$ is injective, and hence $H_2(\mathrm{Aut}^{+}F_n,M)=0$.
Then, considering the homological Lyndon-Hochsild-Serre spectral sequence of (\ref{eq200}), we obtain
$H_2(\mathrm{Aut}\,F_n,M)=0$.

In Section {\rmfamily \ref{2}}, we introduce some tools which we use in our computation in Section {\rmfamily \ref{3}}.
In this paper, a calculation similar to a certain one which we have already mentioned before is often omitted.
(For details, see \cite{Sa2}.)

\tableofcontents

\section{Tools for the computation}\label{2}

In this section, we prepare some tools to compute the second homology groups.
First, we introduce some notation which we use throughout this paper.
Then, we review a finite presentation for $\mathrm{Aut}^{+}F_n$ due to Gersten \cite{Ger} and
reduced finite presentation $\langle X \,|\, R \rangle$ for $\mathrm{Aut}^{+}F_n$ of the Gersten's presentation.
Finally, we show some useful lemmas and equations which are used in Section {\rmfamily \ref{3}}
to reduce the generators of $H_1(\bar{R},M)_{\mathrm{Aut}^{+}F_n}$ where $M= H_L$ and $H_L^*$, and $\bar{R}$ is the normal
closure of $R$ in the free group on the generating set $X$.

\vspace{1em}

Let $F_n$ be a free group of rank $n$ with generators $\{x_1, \ldots , x_n \}$.
In this paper, the group $\mathrm{Aut}\,F_n$ acts on $F_n$ on the right.
For any $\sigma \in \mathrm{Aut}\,F_n$ and $x \in F_n$, the action of $\sigma$ on $x$ is denoted by
$x^{\sigma}$. The elements $x_i^{\pm1} \in F_n$, ($1 \leq i \leq n$), are called letters of $F_n$.
Let $H$ be the abelianization of $F_n$ and $H^*:=\mathrm{Hom}_{\Z}(H,\Z)$ the dual group of $H$.
We remark that although the group $H^*$ is isomorphic to $H$ as a free abelian group, both group are not isomorphic
as an $\mathrm{Aut}\,F_n$-module.
For each generator $x_i \in F_n$, ($1 \leq i \leq n$), set $e_i:=[x_i] \in H$ where $[x]$ means the coset calss of $x$
modulo the commutator subgroup of $F_n$. Then $\{e_1, \ldots , e_n \}$ is a $\Z$-basis of $H$. Let denote
$\{e_1^*, \ldots , e_n^* \}$ the dual basis of it. In general,
in group (co)homology theory, actions of groups on modules are understood to be left actions.
So we consider the modules $H$ and $H^*$ as left $\mathrm{Aut}\,F_n$-modules in a way
$\sigma \cdot x := x^{\sigma^{-1}}$ for $\sigma \in \mathrm{Aut}\,F_n$ and $x \in H$ or $H^*$.

Now, for any letters $a$ and $b$ such that $a \neq b^{\pm1}$,
let $E_{ab}$ be an automorphism of $F_n$ defined by the rule
\[ E_{ab} : \begin{cases}
               a &\mapsto ab, \\
               c &\mapsto c, \hspace{1em} c \neq a^{\pm1}.
              \end{cases}\]
Clearly, we see ${E_{ab}}^{-1}=E_{ab^{-1}}$.
Automorphisms of $F_n$ of $E_{ab}$ type are called Nielsen automorphisms.
In this paper, for simplicity, we write $E_{i^{\pm1}j^{\pm1}}$ for $E_{x_i^{\pm1}x_j^{\pm1}}$.

The actions of $E_{i^{\pm1}j}$ on $e_k$ and $e_k^*$
are given by
\[ E_{i^{\pm1}j} \cdot e_k := [x_k^{{E_{i^{\pm1}j^{-1}}}}] =
                      \begin{cases}
                         e_i \mp e_j, & k=i, \\
                         e_k, & k \neq i
                      \end{cases}\]
and
\[ E_{i^{\pm1}j} \cdot e_k^* =
                      \begin{cases}
                         e_j^* \pm e_i^*, & k=j, \\
                         e_k^*, & k \neq j
                      \end{cases}\]
respectively.
An automorphism $w_{ab} := E_{ba} E_{a^{-1}b} E_{b^{-1}a^{-1}}$ is called a monomial automorphism $a \mapsto b^{-1}$,
$b \mapsto a$. We see ${w_{ab}}^{-1}=w_{ab^{-1}}$, and write
$w_{i^{\pm1}j^{\pm1}}$ for $w_{w_i^{\pm1}x_j^{\pm1}}$.

\vspace{1em}

Let $\rho : \mathrm{Aut}\,F_n \twoheadrightarrow GL(n,\Z)$ be a natural homomorphism induced from the action of $\mathrm{Aut}\,F_n$
on $H$, and $\mathrm{det} : GL(n,\Z) \longrightarrow \{\pm 1\}$ the determinant homomorphism.
The kernel $\mathrm{Aut}^{+}F_n$ of the composition map $\mathrm{det} \circ \rho$
is called the special automorphism group of a free group.
Here we review a finite presentation for $\mathrm{Aut}^{+}F_n$ due to Gersten. He \cite{Ger} showed

\begin{thm}[Gersten \cite{Ger}]\label{LG}
For $n \geq 3$, the group $\mathrm{Aut^{+}}F_n$ has a finite presentation whose generators are $E_{ab}$ subject to relators:
  \begin{description}
    \item[(R1)] ${E_{ab}} E_{a{b^{-1}}}$,
    \item[(R2)] $[E_{ab}, E_{cd}], \,\, \mathrm{for} \,\, a \neq c,d^{\pm1} \,\, \mathrm{and} \,\,  b \neq c^{\pm1}$,
    \item[(R3)] $[E_{ab}, E_{bc}] E_{ac^{-1}}, \,\, \mathrm{for} \,\, a \neq c^{\pm1}$,
    \item[(R4)] $w_{ab} w_{a^{-1} b}$,
    \item[(R5)] ${w_{ab}}^4$.
  \end{description}
\end{thm}
Here $[\,\, , \,\,]$ denotes the commutator bracket defined by $[x,y]:=xyx^{-1}y^{-1}$.
In this paper we often use fundamental formulae of commutators
\begin{equation}
 [x,yz]= [x,y][x,z][[z,x],y], \hspace{1em} [xy,z]=[x,[y,z]][y,z][x,z]. \label{f1}
\end{equation}
We call the relators above the Gersten's relators.
In our paper \cite{Sat}, using Tietze transformations, we reduced the Gersten's presentation to
\begin{lem}\label{L0}
For $n \geq 3$, the group $\mathrm{Aut}^{+}F_n$ has a finite presentation whose
generators are $E_{i^{\pm1}j}$ subject to the relators:
\begin{description}
    \item[(R2-1)] $[E_{ij}, E_{{i}^{-1}j}]$,
    \item[(R2-2)] $[E_{ij}, E_{kj}]$,
    \item[(R2-3)] $[E_{{i}^{-1} j}, E_{kj}]$,
    \item[(R2-4)] $[E_{{i}^{-1}j}, E_{{k}^{-1}j}]$,
    \item[(R2-5)] $[E_{ij}, E_{{i}^{-1}k}]$,
    \item[(R2-6)] $[E_{ij}, E_{kl}]$,
    \item[(R2-7)] $[E_{{i}^{-1}j}, E_{kl}]$,
    \item[(R2-8)] $[E_{{i}^{-1}j}, E_{{k}^{-1}l}]$,
    \item[(R3-1)] $[E_{ik}, E_{kj}] E_{ij^{-1}}$,
    \item[(R3-2)] $[E_{i{k}^{-1}}, E_{{k}^{-1}j}] E_{ij^{-1}}$,
    \item[(R3-3)] $[E_{{i}^{-1}k}, E_{kj}] E_{{i}^{-1}j^{-1}}$,
    \item[(R3-4)] $[E_{{i}^{-1} {k}^{-1}}, E_{{k}^{-1}j}] E_{{i}^{-1}j^{-1}}$,
    \item[(R4-1)] $w_{ij} w_{{i}^{-1}j}$,
    \item[(R5-1)] ${w_{ij}}^4$
  \end{description}
where $i$, $j$, $k$ and $l$ are distinct elements of $\{1, \ldots , n\}$.
\end{lem}
In Section {\rmfamily \ref{3}}, we use this reduced presentation to compute the twisted second homology groups.
In the computation of the second homology groups,

\vspace{0.5em}

Let $X$ and $R$ be the set of generators and relators of the reduced presentation for $\mathrm{Aut}^{+}F_n$ introduced in
Lemma {\rmfamily \ref{L0}} respectively.
In the following, we study relations among the relators of the presentation $\langle X \,|\, R \rangle$,
which is often required in the computation of the second homology groups.
Let $F$ be the free group on $X$, and $\bar{R}$ the normal closure of $R$ in $F$.
Here we define elements $r_{ac}(b)$ and $h_{ab}$ of $F$ to be
\[ r_{ac}(b) := [E_{ab},E_{bc}] \, E_{ac^{-1}} \hspace{1em} \mathrm{for} \hspace{0.5em}
       a \neq b^{\pm1}, c^{\pm1} \hspace{0.5em} \mathrm{and} \hspace{0.5em}  b \neq c^{\pm1} \]
and
\[ h_{ab} := w_{ab} w_{a^{-1}b} \hspace{1em} \mathrm{for} \hspace{0.5em} a \neq b^{\pm1} \]
respectively. Since $r_{ac}(b)$ and $h_{ab}$ are the one of relators of the Gersten's presentation,
we see that these elemets are in $\bar{R}$.
In this paper, we write $r_{i^{\pm1}j^{\pm1}}(k^{\pm1})$ and $h_{i^{\pm}j^{\pm1}}$ for
$r_{x_i^{\pm1}x_j^{\pm1}}(x_k^{\pm1})$ and $h_{x_i^{\pm}x_j^{\pm1}}$ respectively.

For letters $a,b,c$ and $d$, we consider an element $({w_{ab}}^{-1} E_{cd} w_{ab})^{-1} {E_{c^{\sigma}d^{\sigma}}}$ of
$\bar{R}$ where $\sigma$ is the monomial map defined by $w_{ab}$.
More precisely, we study how the elements $({w_{ab}}^{-1} E_{cd} w_{ab})^{-1} {E_{c^{\sigma}d^{\sigma}}}$
are rewritten with the relators of the Gersten's presentation.
First, we consider the case $\sharp \{a^{\pm1},b^{\pm1},$ $ c^{\pm1}, d^{\pm1} \}=6$.
\begin{lem}\label{LA}
For letters $a,b,c,d$, we have
\begin{enumerate}
\item[(i)] if $c=a^{-1}$,
  \[\begin{split}
    ({w_{ab}}^{-1} & E_{a^{-1}d} w_{ab})^{-1} {E_{bd}} \\
        &= E_{b^{-1}a} E_{a^{-1}b^{-1}} \, r_{bd^{-1}}(a^{-1}) \, E_{a^{-1}b} E_{b^{-1}a^{-1}} \\
        & \hspace{2em} \cdot E_{b^{-1}a} E_{bd^{-1}} E_{a^{-1}b^{-1}} \, {r_{a^{-1}d}(b)}^{-1} \,
                                  E_{a^{-1}b} E_{bd} E_{b^{-1}a^{-1}} \\
        & \hspace{2em} \cdot [E_{b^{-1}a},E_{bd^{-1}}].
    \end{split}\]
\item[(ii)] if $c=b^{-1}$,
 \[\begin{split}
    ({w_{ab}}^{-1} & E_{b^{-1}d} w_{ab})^{-1} {E_{a^{-1}d}} \\
        &= E_{b^{-1}a} E_{a^{-1}b^{-1}} \, [E_{ba^{-1}},E_{b^{-1}d^{-1}}] \, E_{a^{-1}b} E_{b^{-1}a^{-1}} \\
        & \hspace{2em} \cdot E_{b^{-1}a} \, {r_{a^{-1}d^{-1}}(b^{-1})} \, E_{b^{-1}a^{-1}} \\
        & \hspace{2em} \cdot E_{a^{-1}d^{-1}} E_{b^{-1}a} \, r_{b^{-1}d}(a^{-1}) \, E_{b^{-1}a^{-1}} E_{a^{-1}d}.
   \end{split}\]
\item[(iii)] if $d=a$,
  \[\begin{split}
    ({w_{ab}}^{-1} & E_{ca} w_{ab})^{-1} {E_{cb^{-1}}} \\
        &= E_{b^{-1}a} E_{a^{-1}b^{-1}} \, [E_{ba^{-1}},E_{ca^{-1}}] \, E_{a^{-1}b} E_{b^{-1}a^{-1}} \\
        & \hspace{2em} \cdot E_{b^{-1}a} E_{cb} \, {r_{cb^{-1}}(a^{-1})}^{-1} \,
                                  E_{cb^{-1}} E_{b^{-1}a^{-1}} \\
        & \hspace{2em} \cdot E_{b^{-1}a} E_{cb} \, {r_{ca^{-1}}(b^{-1})}^{-1} \,
                                  E_{cb^{-1}} E_{b^{-1}a^{-1}}.
    \end{split}\]
\item[(iv)] if $d=a^{-1}$,
  \[\begin{split}
    ({w_{ab}}^{-1} & E_{ca^{-1}} w_{ab})^{-1} {E_{cb}} \\
        &= E_{b^{-1}a} E_{a^{-1}b^{-1}} \, [E_{ba^{-1}},E_{ca}] \, E_{a^{-1}b} E_{b^{-1}a^{-1}} \\
        & \hspace{2em} \cdot E_{b^{-1}a} E_{ca} \, {r_{cb^{-1}}(a^{-1})}^{-1} \,
                                  E_{ca^{-1}} E_{b^{-1}a^{-1}} \\
        & \hspace{2em} \cdot E_{b^{-1}a} E_{ca} \, {r_{ca^{-1}}(b^{-1})}^{-1} \,
                                  E_{ca^{-1}} E_{b^{-1}a^{-1}}.
    \end{split}\]

\item[(v)] if $d=b$,
  \[\begin{split}
    ({w_{ab}}^{-1} & E_{cb} w_{ab})^{-1} {E_{ca}} \\
        &= E_{b^{-1}a} E_{a^{-1}b^{-1}} E_{cb^{-1}} \, r_{ca^{-1}}(b) \, E_{cb} E_{a^{-1}b} E_{b^{-1}a^{-1}} \\
        & \hspace{2em} \cdot E_{b^{-1}a} E_{a^{-1}b^{-1}} E_{cb^{-1}} \, {r_{cb}(a^{-1})} \,
                                  E_{cb} E_{a^{-1}b} E_{b^{-1}a^{-1}} \\
        & \hspace{2em} \cdot [E_{b^{-1}a},E_{ca^{-1}}].
    \end{split}\]

\item[(vi)] if $d=b^{-1}$,
  \[\begin{split}
    ({w_{ab}}^{-1} & E_{cb^{-1}} w_{ab})^{-1} {E_{ca^{-1}}} \\
        &= E_{b^{-1}a} E_{a^{-1}b^{-1}} E_{ca} \, {r_{ca^{-1}}(b)}^{-1} \, E_{ca^{-1}} E_{a^{-1}b} E_{b^{-1}a^{-1}} \\
        & \hspace{2em} \cdot E_{b^{-1}a} E_{ca} \, {r_{cb^{-1}}(a^{-1})} \, E_{ca^{-1}} E_{b^{-1}a^{-1}} \\
        & \hspace{2em} \cdot E_{b^{-1}a} E_{ca} \, [E_{cb^{-1}},E_{a^{-1}b^{-1}}] \, E_{ca^{-1}} E_{b^{-1}a^{-1}} \\
        & \hspace{2em} \cdot [E_{b^{-1}a},E_{ca}].
    \end{split}\]

\end{enumerate}
\end{lem}
Since these equations follows from easy calculations, we omit the details.
In the case where $c=a$ or $c=b$, observing
\begin{equation}
 \begin{split}
    ({w_{ab}}^{-1} E_{cd}  w_{ab})^{-1} {E_{c^{\sigma}d^{\sigma}}}
       & = ({w_{ab}}^{-1} E_{cd^{-1}} h_{ab} E_{cd} w_{ab}) \cdot ({w_{ab}}^{-1} {h_{ab}}^{-1} w_{ab}) \\
       & \hspace{2em} \cdot ({w_{a^{-1}b^{-1}}}^{-1} E_{cd} w_{a^{-1}b^{-1}})^{-1} {E_{c^{\sigma}d^{\sigma}}},
  \end{split} \label{eq2}
\end{equation}
and Lemma {\rmfamily \ref{LA}} above, we see that
$({w_{ab}}^{-1} E_{cd}  w_{ab})^{-1} {E_{c^{\sigma}d^{\sigma}}}$ is also rewritten with the relators of the Gersten's presentation.
For the case $\sharp \{a^{\pm1},b^{\pm1},$ $ c^{\pm1}, d^{\pm1} \}=8$, we have
\begin{lem}\label{LB}
\[\begin{split}
    ({w_{ab}}^{-1} E_{cd} w_{ab})^{-1} {E_{cd}}
        &= E_{b^{-1}a} E_{a^{-1}b^{-1}} \, [E_{ba^{-1}},E_{cd^{-1}}] \, E_{a^{-1}b} E_{b^{-1}a^{-1}} \\
        & \hspace{2em} \cdot E_{b^{-1}a} \, [E_{a^{-1}b^{-1}},E_{cd^{-1}}] \, E_{b^{-1}a^{-1}}
        \cdot [E_{b^{-1}a},E_{cd^{-1}}].
    \end{split}\]
\end{lem}

Next, we consider how rewrite the relators $[E_{ab},E_{cd}]$, $r_{ac}(b)$ and $h_{ab}$ of the Gersten's presentation with
the relators (R2-1), $\ldots$ , (R4-1) of the reduced presentation.
First, by an easy calculation, we see that the (R2) type relator $[E_{ab},E_{cd}]$ is rewritten as a conjugate of one of the relators
(R2-1), $\ldots$ , (R2-8). For example,
\[ [E_{ij^{-1}},E_{kj^{-1}}] = E_{ij^{-1}} E_{kj^{-1}} \, [E_{ij},E_{kj}] \, E_{ij} E_{kj}. \]
For the relators $r_{ac}(b)$ and $h_{ab}$, we use
\begin{lem}\label{LC}
\hspace*{\fill}\ 

\begin{enumerate}
\item[(i)] $r_{ac^{-1}}(b) = (E_{bc^{-1}}E_{ac^{-1}} {r_{ac}(b)}^{-1} E_{ac} E_{bc}) \cdot [E_{bc^{-1}},E_{ac^{-1}}]$,
\item[(ii)] $h_{a^{-1}b} = {w_{ab}}^{-1} h_{ab} w_{ab}$, \hspace{1em} $h_{ab^{-1}} = {w_{ab}}^{-1} {h_{ab}}^{-1} w_{ab}$.
\end{enumerate}
\end{lem}

\vspace{0.5em}

Finally, we consider two type of equations induced from elements of $\bar{R}$:
\[ ({w_{ab}}^{-1} r_{cd}(e) w_{ab})^{-1} r_{c^{\sigma} {d^{\sigma}}}(e^{\sigma}) \hspace{1em} \mathrm{and} \hspace{1em}
   ({w_{ab}}^{-1} h_{cd} w_{ab})^{-1} h_{c^{\sigma} {d^{\sigma}}} \]
where $\sigma$ is the monomial map defined by $w_{ab}$.
Observe
\begin{equation}
\begin{split}
   r_{c^{\sigma} {d^{\sigma}}}(e^{\sigma}) & =
     s_1^{-1} \cdot ({w_{ab}}^{-1} E_{ce} w_{ab} \,\, s_2^{-1} \,\, {w_{ab}}^{-1} {E_{ce}}^{-1} w_{ab}) \\
   & \hspace{2em} \cdot ({w_{ab}}^{-1} r_{cd}(e) w_{ab}) \\
   & \hspace{2em} \cdot ({w_{ab}}^{-1} E_{cd} E_{ed} w_{ab} \,\, s_1 \,\, {w_{ab}}^{-1} {E_{ed}}^{-1} {E_{cd}}^{-1} w_{ab}) \\
   & \hspace{2em} \cdot ({w_{ab}}^{-1} E_{cd} w_{ab} \,\, s_2 \,\, {w_{ab}}^{-1} {E_{cd}}^{-1} w_{ab}) \cdot s_3
\end{split} \label{eq21}
\end{equation}
where
\[\begin{split}
   s_1 & := ({w_{ab}}^{-1} E_{ce^{-1}} w_{ab})^{-1} {E_{c^{\sigma} {e^{\sigma}}}}^{-1}, \hspace{1em}
   s_2 := ({w_{ab}}^{-1} E_{ed^{-1}} w_{ab})^{-1} {E_{e^{\sigma} {d^{\sigma}}}}^{-1}, \\
   s_3 & := ({w_{ab}}^{-1} E_{cd^{-1}} w_{ab})^{-1} {E_{c^{\sigma} {d^{\sigma}}}}^{-1}. \\
  \end{split}\]
Considering the equation (\ref{eq21}) in $\bar{R}^{\mathrm{ab}}$,
and tensoring it with $e_p$ in $\bar{R}^{\mathrm{ab}} \otimes_{\Z} H_L$, we obtain
\[\begin{split}
   r_{c^{\sigma} {d^{\sigma}}}(e^{\sigma}) \otimes e_p & =
     s_1^{-1} \otimes e_p + ({w_{ab}}^{-1} E_{ce} w_{ab} \,\, s_2^{-1} \,\, {w_{ab}}^{-1} {E_{ce}}^{-1} w_{ab}) \otimes e_p \\
   & \hspace{2em} + ({w_{ab}}^{-1} r_{cd}(e) w_{ab}) \otimes e_p \\
   & \hspace{2em} + ({w_{ab}}^{-1} E_{cd} E_{ed} w_{ab} \,\, s_1 \,\, {w_{ab}}^{-1} {E_{ed}}^{-1} {E_{cd}}^{-1} w_{ab}) \otimes e_p \\
   & \hspace{2em} + ({w_{ab}}^{-1} E_{cd} w_{ab} \,\, s_2 \,\, {w_{ab}}^{-1} {E_{cd}}^{-1} w_{ab}) \otimes e_p + s_3 \otimes e_p
\end{split}\]
For convenience, we denote this equation by $(a,b,c,d,e) \otimes e_p$. Similarly, we define $(a,b,c,d,e) \otimes e_p^*$.

We also consider
\begin{equation}
\begin{split}
   h_{c^{\sigma} {d^{\sigma}}} & =
     t_3 \cdot ({w_{ab}}^{-1} E_{dc} w_{ab} \,\, t_2 \,\, {w_{ab}}^{-1} {E_{dc}}^{-1} w_{ab}) \\
   & \hspace{2em} \cdot ({w_{ab}}^{-1} E_{dc} E_{c^{-1}d} w_{ab} \,\, t_1 \,\, {w_{ab}}^{-1} {E_{c^{-1}d}}^{-1} {E_{dc}}^{-1} w_{ab}) \\
   & \hspace{2em} \cdot ({w_{ab}}^{-1} {E_{d^{-1}c}}^{-1} {E_{cd}}^{-1} w_{ab}
                             \,\, t_4 \,\, {w_{ab}}^{-1} E_{cd} {E_{d^{-1}c}} w_{ab}) \\
   & \hspace{2em} \cdot ({w_{ab}}^{-1} {E_{d^{-1}c}}^{-1} w_{ab}
                             \,\, t_5 \,\, {w_{ab}}^{-1} {E_{d^{-1}c}} w_{ab}) \cdot t_6
\end{split} \label{eq41}
\end{equation}
where
\[\begin{split}
   t_1 & :=  E_{(d^{-1})^{\sigma} (c^{-1})^{\sigma}} (w_{ab}^{-1} E_{d^{-1}c^{-1}} w_{ab})^{-1},
                  \hspace{2.8em} t_2 := E_{(c^{-1})^{\sigma} d^{\sigma}} (w_{ab}^{-1} E_{c^{-1}d} w_{ab})^{-1}, \\
   t_3 & := E_{d^{\sigma} c^{\sigma}} (w_{ab}^{-1} E_{dc} w_{ab})^{-1},
                        \hspace{7.25em} t_4 := (w_{ab}^{-1} E_{dc^{-1}} w_{ab})^{-1} E_{d^{\sigma} (c^{-1})^{\sigma}} , \\
   t_5 & := (w_{ab}^{-1} E_{cd} w_{ab})^{-1} E_{c^{\sigma} d^{\sigma}},
                        \hspace{7.25em} t_6 := (w_{ab}^{-1} E_{d^{-1}c} w_{ab})^{-1} E_{(d^{-1})^{\sigma} c^{\sigma}}.
  \end{split}\]
For convenience, we denote by $\{ a,b,c,d \} \otimes e_p^*$
the equations obtained by considering (\ref{eq41}) in $\bar{R}^{\mathrm{ab}}$, and tensoring it with
$e_p^*$ in $\bar{R}^{\mathrm{ab}} \otimes_{\Z} H_L^*$.
We often use these equations in Section {\rmfamily \ref{3}} to reduce the generators of
$H_1(\bar{R},M)_{\mathrm{Aut}^{+}F_n}$ for $M= H_L$ and $H_L^*$.

\vspace{1em}

\section{The Proof of the main theorem}\label{3}

First we consider $H_2(\mathrm{Aut}^{+}F_n,H_L)$ for $n \geq 6$.
Let $F$, $R$ and $\bar{R}$ be as above.
Then we have an exact sequence
\[ 1 \rightarrow \bar{R} \rightarrow F \rightarrow \mathrm{Aut}^{+}F_n \rightarrow 1. \]
This sequence induces a homological five-term exact sequence
\[\begin{split}
    H_2(F,H_L) \rightarrow H_2(\mathrm{Aut}^{+}F_n, & H_L) \rightarrow H_1(\bar{R},H_L)_{\mathrm{Aut}^{+}F_n} \\
                     & \rightarrow H_1(F,H_L) \rightarrow H_1(\mathrm{Aut}^{+}F_n,H_L) \rightarrow 0.
  \end{split}\]
of $\Z$-modules. 
Since a $\Z$-equivariant homomorphism between $L$-modules is naturally considered as a $L$-equivariant homomorphism, we see that
this sequence is an $L$-equivariant exact sequence.
Since $F$ is a free group, 
$H_2(F,H_L)=0$. Furthermore $H_1(\mathrm{Aut}^{+}F_n,H_L)=0$ from our results of \cite{Sat}.
Hence we have an $L$-equivariant short exact sequence
\[ 0 \rightarrow H_2(\mathrm{Aut}^{+}F_n, H_L) \rightarrow H_1(\bar{R},H_L)_{\mathrm{Aut}^{+}F_n}
      \rightarrow H_1(F,H_L) \rightarrow 0.\]
Since $F$ is a free group of rank $2(n^2-n)$, and since $H_L$ is a free $L$-module of rank $n$,
we see $H^1(F,H_L)$ is a free $L$-module of rank $2n(n^2-n)-n$. Hence, by universal coefficients theorem, we have
$H_1(F,H_L) \simeq L^{\oplus \{2n(n^2-n)-n \}}$.
Since $L$ is a principal ideal domain, we can apply the structure theorem to any finitely generated $L$-modules.
Therefore our required result $H_2(\mathrm{Aut}^{+}F_n,H_L)=0$ follows from
\begin{prop}\label{P1}
For $n \geq 6$,
\[ H_1(\bar{R},H_L)_{\mathrm{Aut}^{+}F_n} \simeq L^{\oplus \{2n(n^2-n)-n \} }. \]
\end{prop}
We prove this proposition in Subsection {\rmfamily \ref{H}}.
Then, observing the homological Lyndon-Hochschild-Serre spectral sequence of
\begin{equation}
 1 \rightarrow \mathrm{Aut^{+}}F_n \rightarrow \mathrm{Aut}\,F_n \rightarrow \{\pm1\} \rightarrow 1, \label{ex1}
\end{equation}
we obtain $H_2(\mathrm{Aut}\,F_n,H_L)=0$ for $n \geq 6$.

\vspace{1em}

Next we consider $H_2(\mathrm{Aut}^{+}F_n,H_L^*)$ for $n \geq 6$.
Similarly, we obtain a homological five-term exact sequence
\[\begin{split}
   H_2(F,H_L^*) \rightarrow H_2(\mathrm{Aut}^{+}F_n, & H_L^*) \rightarrow H_1(\bar{R},H_L^*)_{\mathrm{Aut}^{+}F_n} \\
                     & \rightarrow H_1(F,H_L^*) \rightarrow H_1(\mathrm{Aut}^{+}F_n,H_L^*) \rightarrow 0,
  \end{split}\]
of $L$-modules, and from the results of \cite{Sat},
\[ 0 \rightarrow H_2(\mathrm{Aut}^{+}F_n, H_L^*) \rightarrow H_1(\bar{R},H_L^*)_{\mathrm{Aut}^{+}F_n}
      \rightarrow H_1(F,H_L^*) \rightarrow L \rightarrow 0.\]
Since we have $H_1(F,H_L^*) \simeq L^{\oplus \{2n(n^2-n)-n \}}$, our required result $H_2(\mathrm{Aut}^{+}F_n,H_L^*)=0$ follows from
\begin{prop}\label{P2}
For $n \geq 6$,
\[ H_1(\bar{R},H_L^*)_{\mathrm{Aut}^{+}F_n} \simeq L^{\oplus \{2n(n^2-n)-n-1 \} }. \]
\end{prop}
We prove this proposition in Subsection {\rmfamily \ref{H^*}}.
Then observing the homological Lyndon-Hochschild-Serre spectral sequence of (\ref{ex1}),
we obtain $H_2(\mathrm{Aut}\,F_n,H_L^*)=0$ for $n \geq 6$.

\subsection{The proof of Proposotion {\rmfamily \ref{P1}}}\label{H}
\hspace*{\fill}\ 

\vspace{0.5em}

In this subsection, we prove Proposotion {\rmfamily \ref{P1}}.
Since the map $H_1(\bar{R},H_L)_{\mathrm{Aut}^{+}F_n} \rightarrow H_1(F,H_L) = L^{\oplus \{2n(n^2-n)-n \}}$
is surjective, $H_1(\bar{R},H_L)_{\mathrm{Aut}^{+}F_n}$ contains a free $L$-submodule which rank is
greater than or equal to $2n(n^2-n)-n$.
To show it is just $2n(n^2-n)-n$, it suffices to show that
$H_1(\bar{R},H_L)_{\mathrm{Aut}^{+}F_n}$ is generated by just $2n(n^2-n)-n$ elements.
Let $\bar{R}^{\mathrm{ab}}$ be the abelianization of $\bar{R}$. We also denote by $r$ the coset class of $r \in \bar{R}$.
By definifition, we have $H_1(\bar{R},H_L)_{\mathrm{Aut}^{+}F_n} = \bar{R}^{\mathrm{ab}} {\otimes}_{\mathrm{Aut}^{+}F_n} H_L$,
and see that
\[ \mathfrak{E} := \{ r \otimes e_p \,|\, r \in R, \,\, 1 \leq p \leq n \} \]
is a generating set of $\bar{R}^{\mathrm{ab}} {\otimes}_{\mathrm{Aut}^{+}F_n} H_L$ as a $L$-module.
In the following, we reduce the elements of $\mathfrak{E}$.
We use $\equiv$ for the equality in $\bar{R}^{\mathrm{ab}} {\otimes}_{\mathrm{Aut}^{+}F_n} H_L$.

\vspace{1em}

{\bf Step 0.}
In the reduction of the generators of $\bar{R}^{\mathrm{ab}} {\otimes}_{\mathrm{Aut}^{+}F_n} H_L$, we often use the following
lemmas.

\begin{lem}\label{L1}
For $n \geq 3$,
\[\begin{split}
   (E_{{i}^{\pm 1}j} \, r \, {E_{{i}^{\pm 1}j^{-1}}}) \otimes e_p & \equiv
              \begin{cases}
                r \otimes e_p, & \hspace{1em} p \neq i, \\
                r \otimes e_i \pm r \otimes e_j, & \hspace{1em} p=i,
              \end{cases} \\
   ({E_{{i}^{\pm 1}j^{-1}}} \, r \, E_{{i}^{\pm 1}j}) \otimes e_p & \equiv
              \begin{cases}
                r \otimes e_p, & \hspace{1em} p \neq i, \\
                r \otimes e_i \mp r \otimes e_j, & \hspace{1em} p=i.
              \end{cases} \\
  \end{split}\]
\end{lem}

\textit{Proof of Lemma {\rmfamily \ref{L1}}}
For any $\sigma \in \mathrm{Aut}^{+}F_n$, $r \in \bar{R}^{\mathrm{ab}}$ and $h \in H$,
we have $r \cdot \sigma \otimes h \equiv r \otimes \sigma \cdot h$.
Then observing the equation $r \cdot \sigma \otimes h = \sigma^{-1} r \sigma \otimes h$ induced from the definition of the action
of $\mathrm{Aut}^{+}F_n$ on $\bar{R}^{\mathrm{ab}}$,
we obtain the required results by substituting $\sigma = E_{i^{\pm1} j^{\pm1}}$ and $h=e_p$ . $\square$

\begin{cor}\label{C1}
For $n \geq 3$,
\[\begin{split}
   [E_{{i}^{\pm1} j}, r] \otimes e_p & \equiv
              \begin{cases}
                0, & \hspace{1em} p \neq i, \\
                \pm r \otimes e_j, & \hspace{1em} p=i,
              \end{cases} \\
   [{E_{{i}^{\pm1} j^{-1}}}, r] \otimes e_p & \equiv 
              \begin{cases}
                0, & \hspace{1em} p \neq i, \\
                \mp r \otimes e_j, & \hspace{1em} p=i
              \end{cases}
  \end{split}\]
\end{cor}

\textit{Proof of Corollary {\rmfamily \ref{C1}}} 
Observing
\[ [\sigma,r] \otimes e_p = \sigma r \sigma^{-1} r^{-1} \otimes e_p = (\sigma r \sigma^{-1}) \otimes e_p - r \otimes e_p \]
for $\sigma \in F$ and $r \in \bar{R}^{\mathrm{ab}}$, and Lemma {\rmfamily \ref{L1}},
we immediately obtain the required results. $\square$

\vspace{0.5em}

Similarly, we have
\begin{lem}\label{L2}
For $n \geq 3$,
\[(w_{{i}^{\pm 1}j} \, r \, {w_{{i}^{\pm 1}j}}^{-1}) \otimes e_p \equiv
              \begin{cases}
                r \otimes e_p, & \hspace{1em} p \neq i,j, \\
                \mp r \otimes e_j, & \hspace{1em} p=i, \\
                \pm r \otimes e_i, & \hspace{1em} p=j,
              \end{cases}\]
\[(w_{{i}^{\pm 1}j^{-1}} \, r \, {w_{{i}^{\pm 1}j}}) \otimes e_p \equiv
              \begin{cases}
                r \otimes e_p, & \hspace{1em} p \neq i,j, \\
                \pm r \otimes e_j, & \hspace{1em} p=i, \\
                \mp r \otimes e_i, & \hspace{1em} p=j.
              \end{cases}\]
\end{lem}
\begin{cor}\label{C2}
For $n \geq 3$,
\[ [w_{{i}^{\pm1}j}, r] \otimes e_p \equiv
   \begin{cases}
     0,  \hspace{1em} & p \neq i,j, \\
     -r \otimes e_i \mp r \otimes e_j, \hspace{1em} & p = i, \\
     \pm r \otimes e_i -r \otimes e_j, \hspace{1em} & p = j. \\
   \end{cases}\]
\[ [{w_{{i}^{\pm1}j}}^{-1}, r] \otimes e_p \equiv
   \begin{cases}
     0,  \hspace{1em} & p \neq i,j, \\
     -r \otimes e_i \pm r \otimes e_j, \hspace{1em} & p = i, \\
     \mp r \otimes e_i -r \otimes e_j, \hspace{1em} & p = j. \\
   \end{cases}\]
\end{cor}

\vspace{0.5em}

Considering any relator of (R2) of the Gersten's presentation is conjugate to one of the relator of (R2-1), $\ldots$ , (R2-8),
or considering Lemma {\rmfamily \ref{LC}}, for any relator $r=$ (R2), (R3) and (R4),
we can rewrite a element $r \otimes e_p$ with the relators (R2-1), $\ldots$ , (R4-1)
using Lemmas {\rmfamily \ref{L1}} and {\rmfamily \ref{L2}}.
The computation is easiest explained with examples, so we give three examples.
\[\begin{split}
   [E_{ij^{-1}},E_{kj^{-1}}] \otimes e_p & =( E_{ij^{-1}} E_{kj^{-1}} \, [E_{ij},E_{kj}] \, E_{ij} E_{kj} ) \otimes e_p, \\
  & \equiv  \begin{cases}
               [E_{ij},E_{kj}] \otimes e_p,  \hspace{1em} & p \neq i, k, \\
               [E_{ij},E_{kj}] \otimes e_p - [E_{ij},E_{kj}] \otimes e_j, \hspace{1em} & p = i,k. \\
            \end{cases}
  \end{split}\]

\[\begin{split}
   r_{ij^{-1}}(k) \otimes e_p & =
              \{ (E_{kj^{-1}}E_{ij^{-1}} {r_{ij}(k)}^{-1} E_{ij} E_{kj}) \cdot [E_{kj^{-1}},E_{ij^{-1}}] \} \otimes e_p, \\
    & \equiv  \begin{cases}
               - r_{ij}(k) \otimes e_p + [E_{ij},E_{kj}] \otimes e_p,  \hspace{1em} & p \neq i, k, \\
               - r_{ij}(k) \otimes (e_p-e_j) + [E_{ij},E_{kj}] \otimes (e_p - e_j), \hspace{1em} & p = i,k. \\
            \end{cases}
  \end{split}\]
\[\begin{split}
   h_{i^{-1}j} \otimes e_p & = ({w_{ij}}^{-1} h_{ij} w_{ij}) \otimes e_p
      \equiv  \begin{cases}
               h_{ij} \otimes e_p,  \hspace{1em} & p \neq i,j, \\
               h_{ij} \otimes e_j,  \hspace{1em} & p = i, \\
               - h_{ij} \otimes e_i,  \hspace{1em} & p = j.
            \end{cases}
  \end{split}\]

\vspace{1em}

{\bf Step 1.}
First we consider the generators ${w_{ij}}^4 \otimes e_p$.
Observing the Gersten's computation in \cite{Ger}, we see that for any $a,b,c,$ and $d$, the element
$({w_{ab}}^{-1} E_{cd}$ $ w_{ab})^{-1}{E_{c^{\sigma}d^{\sigma}}} \in \bar{R}$
is in the normal closure of (R2-1), $\ldots$ , (R4-1) in $F$, Hence
\[ {w_{ij}}^8 =
   \{ ({w_{ij}}^4 {w_{jk}}^{2} {w_{ij}}^{-4})^{-1} {w_{jk}}^2 \} 
   \cdot \{ ({w_{jk}}^{-2} {w_{ij}}^{-4} {w_{jk}}^{2})^{-1} {w_{ij}}^4) \} \]
is also in it, and we see that
\[ {w_{ij}}^4 \otimes e_p= \frac{1}{2} \, (2 \, {w_{ij}}^4 \otimes e_p) \equiv \frac{1}{2} \, ({w_{ij}}^8 \otimes e_p) \]
is rewritten as a sum of the generator $r \otimes e_p$ for $r=$ (R2-1), $\ldots$ , (R4-1).
Therefore we can remove the generators ${w_{ij}}^4 \otimes e_p$ from the generationg set $\mathfrak{E}$.

\vspace{1em}

{\bf Step 2.}
Here we show that the generators $r \otimes e_p$ for $r=$ (R2-1), $\ldots$ , (R2-8) is zero or equal to one of 
the generators $r_{{i}^{\pm1}j}(k^{\pm1}) \otimes e_p$. We have
\begin{lem}\label{L3}
For $n \geq 6$ and distinct $i,j,k,l$ and $m$,
\begin{enumerate}
\item[(i)] {\bf{(R2-2):}}
\[ [E_{ij}, E_{kj}] \otimes e_p \equiv
   \begin{cases}
     0,  \hspace{1em} & p \neq i,k, \\
     -r_{kj}(l) \otimes e_j, \hspace{1em} & p = i, \\
     r_{ij}(l) \otimes e_j, \hspace{1em} & p = k. \\
   \end{cases}\]
\item[(ii)] {\bf{(R2-3), (R2-4):}}
\[ [E_{{i}^{-1}j}, E_{{k}^{\pm1}j}] \otimes e_p \equiv
   \begin{cases}
     0,  \hspace{1em} & p \neq i,k, \\
     r_{{k}^{\pm1}j}(l) \otimes e_j, \hspace{1em} & p = i, \\
     \mp r_{{i}^{-1}j}(l) \otimes e_j, \hspace{1em} & p = k. \\
   \end{cases}\]
\item[(iii)] {\bf{(R2-1):}}
\[ [E_{ij}, E_{{i}^{-1}j}] \otimes e_p \equiv
   \begin{cases}
     0,  \hspace{1em} & p \neq i, \\
     - r_{ij}(k) \otimes e_j - r_{{i}^{-1}j}(l) \otimes e_j, \hspace{1em} & p = i.
   \end{cases}\]
\item[(iv)] {\bf{(R2-6):}}
\[ [E_{ij}, E_{kl}] \otimes e_p \equiv
   \begin{cases}
     0,  \hspace{1em} & p \neq i,k, \\
     - r_{kl}(m) \otimes e_j, \hspace{1em} & p = i, \\
     r_{ij}(m) \otimes e_l, \hspace{1em} & p = k. \\
   \end{cases}\]
\item[(v)] {\bf{(R2-7), (R2-8):}}
\[ [E_{{i}^{-1}j}, E_{{k}^{\pm1}l}] \otimes e_p \equiv
   \begin{cases}
     0,  \hspace{1em} & p \neq i,k, \\
     r_{{k}^{\pm1}l}(m) \otimes e_j, \hspace{1em} & p = i, \\
     \pm r_{{i}^{-1}j}(m) \otimes e_l, \hspace{1em} & p = k. \\
   \end{cases}\]
\item[(vi)] {\bf{(R2-5):}}
\[ [E_{ij}, E_{{i}^{-1}k}] \otimes e_p \equiv
   \begin{cases}
     0,  \hspace{1em} & p \neq i, \\
     - r_{ij}(l) \otimes e_k - r_{{i}^{-1}k}(m) \otimes e_j, \hspace{1em} & p = i.
   \end{cases}\]
\end{enumerate}
\end{lem}

\textit{Proof of Lemma {\rmfamily \ref{L3}}}
Here we prove (i).
First we consider the case $p \neq i,k$. Since $n \geq 5$, we can choose a number $l \in \{1, \ldots , n \}$ such that
$l \neq i,j,k,p$. Set $r:= E_{ij^{-1}} [E_{il},E_{lj}] \in \bar{R}$. Since $[r,E_{kj}]$ is in $\bar{R}$, and since
$p \neq i,k$, we have
\[ [E_{ij},[r,E_{kj}]] \otimes e_p \equiv 0, \hspace{1em} [r,E_{kj}] \otimes e_p \equiv 0 \]
by Corollary {\rmfamily \ref{C1}}.
Then, using the formula (\ref{f1}) repeatedly, we see
\[\begin{split}
   [E_{ij}, E_{kj}] \otimes e_p
     & \equiv ([E_{ij},[r,E_{kj}]] [r,E_{kj}] [E_{ij}, E_{kj}]) \otimes e_p, \\
     & = [[E_{il},E_{lj}],E_{kj}] \otimes e_p, \\
     & = [ E_{il} E_{lj} {E_{il}}^{-1} {E_{lj}}^{-1}, E_{kj}] \otimes e_p, \\
     & = ([E_{il}, [E_{lj} {E_{il}}^{-1} {E_{lj}}^{-1}, E_{kj}]] 
           [E_{lj} {E_{il}}^{-1} {E_{lj}}^{-1}, E_{kj}] [E_{il}, E_{kj}]) \otimes e_p, \\
     & \equiv ([E_{lj} {E_{il}}^{-1} {E_{lj}}^{-1}, E_{kj}] [E_{il}, E_{kj}]) \otimes e_p, \\
     & \equiv \cdots \cdots, \\
     & \equiv ([{E_{lj}}^{-1}, E_{kj}] [{E_{il}}^{-1}, E_{kj}] [E_{lj}, E_{kj}] [E_{il}, E_{kj}]) \otimes e_p, \\
     & \equiv [{E_{lj}}^{-1}, E_{kj}] \otimes e_p + [{E_{il}}^{-1}, E_{kj}] \otimes e_p + [E_{lj}, E_{kj}] \otimes e_p
                    + [E_{il}, E_{kj}] \otimes e_p.
\end{split}\]
On the other hand, by (\ref{f1}) we have
\[ 1 = [E_{lj} {E_{lj}}^{-1}, E_{kj}] = [E_{lj},[{E_{lj}}^{-1}, E_{kj}]] [{E_{lj}}^{-1}, E_{kj}] [E_{lj}, E_{kj}]. \]
Since $p \neq l$, we see $[E_{lj},[{E_{lj}}^{-1}, E_{kj}]] \otimes e_p \equiv 0$ by Corollary {\rmfamily \ref{C1}},
and hence
\[ [E_{lj}, E_{kj}] \otimes e_p + [{E_{lj}}^{-1}, E_{kj}] \otimes e_p \equiv 0. \]
Similarly, since $p \neq i$,
\[ [E_{il}, E_{kj}] \otimes e_p + [{E_{il}}^{-1}, E_{kj}] \otimes e_p \equiv 0. \]
Therefore we obtain $[E_{ij}, E_{kj}] \otimes e_p \equiv 0.$

\vspace{0.5em}

Next we consider the case $p=i$. Since $[[E_{kj},E_{ij}], r_{kj}(l)] =0$ in ${\bar{R}}^{\mathrm{ab}}$, and since
$[E_{ij}, r_{kj}(l)] \otimes e_i \equiv r_{kj}(l) \otimes e_j$ by Corollary {\rmfamily \ref{C1}}, we have
\[\begin{split}
   [E_{ij},  E_{kj}] \otimes & e_i + r_{kj}(l) \otimes e_j \\
     & \equiv ([E_{ij}, r_{kj}(l)] [E_{ij},  E_{kj}] [[E_{kj},E_{ij}], r_{kj}(l)]) \otimes e_i, \\
     & = [ E_{ij}, [E_{kl},E_{lj}]] \otimes e_i, \\
     & = ([E_{ij},E_{kl}] [E_{ij}, E_{lj} {E_{kl^{-1}}} {E_{lj^{-1}}}]
         [[E_{lj} {E_{kl^{-1}}} {E_{lj^{-1}}}, E_{ij}], E_{kl}]) \otimes e_i, \\
     & \equiv ([E_{ij},E_{kl}] [E_{ij}, E_{lj} {E_{kl^{-1}}} {E_{lj^{-1}}}]) \otimes e_i, \\ 
     & \equiv \cdots \cdots, \\
     & \equiv ([E_{ij},E_{kl}] [E_{ij},{E_{kl^{-1}}}] [E_{ij},E_{lj}] [E_{ij},{E_{lj^{-1}}}]) \otimes e_i \\
     & \equiv ([E_{ij},E_{kl}] + [E_{ij},{E_{kl^{-1}}}]) \otimes e_i + ([E_{ij},E_{lj}] + [E_{ij},{E_{lj^{-1}}}]) \otimes e_i.
  \end{split}\]
Since we have
\[ 1 = [E_{ij},E_{kl} E_{kl^{-1}}] = [E_{ij},E_{kl}] [E_{ij},E_{kl^{-1}}] [[E_{kl^{-1}},E_{ij}], E_{kl}], \]
and since $[[E_{kl^{-1}},E_{ij}], E_{kl}] \otimes e_i \equiv 0$ by Corollary {\rmfamily \ref{C1}}, we see
\[ ([E_{ij},E_{kl}] + [E_{ij},{E_{kl^{-1}}}]) \otimes e_i \equiv 0 .\]
Similarly,
\[ ([E_{ij},E_{lj}] + [E_{ij},{E_{lj^{-1}}}]) \otimes e_i \equiv 0. \]
Hence we obtain $[E_{ij},  E_{kj}] \otimes e_i \equiv - r_{kj}(l) \otimes e_j$.
Furthermore changing the role of $i$ and $k$ in the equation $[E_{kj},  E_{ij}] \otimes e_i \equiv r_{kj}(l) \otimes e_j$,
we also obtain $[E_{ij},  E_{kj}] \otimes e_k \equiv r_{ij}(l) \otimes e_j$. 

\vspace{0.5em}

Similarly, we can show (ii), (iv) and (v). We remark that to show (iv) and (v), we need $n \geq 6$ since we use six distinct
generators of the free group $F_n$.
Then using these results, we obtain (iii) and (vi). Since the calculations are
similar to that above, we leave it to the reader for exercise. (For details, see \cite{Sa2}.)
$\square$

\vspace{0.5em}

By the lemma above, we can remove the generators $r \otimes e_p$ for $r=$ (R2-1), $\ldots$ , (R2-8)
from the generationg set $\mathfrak{E}$.

\vspace{1em}

{\bf Step 3.} Here we consider the generators $r_{{i}^{\pm1}j}(k^{\pm1}) \otimes e_p$ for $p \neq i$.

\vspace{0.5em}

{\bf (3-a)} {\bf The case $p \neq i,k$}.

\vspace{0.5em}

First we consider the case $p=j$.
Observing (i) of Lemma {\rmfamily \ref{L3}}, we see that $r_{ij}(l) \otimes e_j$ doesn't depend on
the choice of a number $l$ such that $l \neq i,j,k$.
On the other hand, since $n \geq 5$, there exists another number $m$ such that $m \neq i,j,k,l$.
Similarly, we have
$[E_{ij}, E_{mj}] \otimes e_m \equiv r_{ij}(k) \otimes e_j \equiv r_{ij}(l) \otimes e_j$ from (i) of Lemma {\rmfamily \ref{L3}}.
This shows that $r_{ij}(l) \otimes e_j$ doesn't depend on the choice of a number $l$ such that $l \neq i,j$.
Futhermore, using the relator $r_{kj}(l^{-1})$ instead of $r_{kj}(l)$ in the proof of (i) of Lemma {\rmfamily \ref{L3}},
we also obtain
\[ [E_{ij}, E_{kj}] \otimes e_p \equiv
   \begin{cases}
     0,  \hspace{1em} & p \neq i,k, \\
     - r_{kj}(l^{-1}) \otimes e_j, \hspace{1em} & p = i, \\
     r_{ij}(l^{-1}) \otimes e_j, \hspace{1em} & p = k. \\
   \end{cases}\]
Hence we can set
\[ r_{ij}(\cdot) \otimes e_j :\equiv  r_{ij}(k) \otimes e_j \equiv r_{ij}(k^{-1}) \otimes e_j \]
for distinct $i$ and $j$.
Similarly, observing (ii) of Lemma {\rmfamily \ref{L3}}, we can set
\[ r_{{i}^{-1}j}(\cdot) \otimes e_j :\equiv  r_{{i}^{-1}j}(k) \otimes e_j \equiv r_{{i}^{-1}j}(k^{-1}) \otimes e_j \]
for distinct $i$ and $j$.

\vspace{0.5em}

For the case $p \neq j$, observing (iv) and (v) of Lemma {\rmfamily \ref{L3}}, we can set
\[\begin{split}
   r_{ij}(\cdot) \otimes e_p :\equiv & r_{ij}(k) \otimes e_p \equiv r_{ij}(k^{-1}) \otimes e_p, \\
   r_{{i}^{-1}j}(\cdot) \otimes e_p :\equiv & r_{{i}^{-1}j}(k) \otimes e_p \equiv r_{{i}^{-1}j}(k^{-1}) \otimes e_p.
  \end{split}\]

\vspace{1em}

{\bf (3-b)} {\bf The case $p=k$.}

\vspace{0.5em}

Set
\[ S_{ijk} := r_{ij}(k) \otimes e_k - r_{ij}(\cdot) \otimes e_k - r_{ik}(\cdot) \otimes e_j. \]
We show that $S_{ijk} \equiv 0$ in $\bar{R}^{\mathrm{ab}} {\otimes}_{\mathrm{Aut}^{+}F_n} H_L$.
For distinct $i,j,k$ and $l$, the equation $(x_l,x_j,x_i,x_j,x_k) \otimes e_k$ is given by
\begin{equation}
\begin{split}
   r_{il}(k) \otimes e_k & =
     s_1^{-1} \otimes e_k + ({w_{lj}}^{-1} E_{ik} w_{lj} \,\, s_2^{-1} \,\, {w_{lj}}^{-1} {E_{ik}}^{-1} w_{lj}) \otimes e_k \\
   & \hspace{2em} + ({w_{lj}}^{-1} r_{ij}(k) w_{lj}) \otimes e_k \\
   & \hspace{2em} + ({w_{lj}}^{-1} E_{ij} E_{kj} w_{lj} \,\, s_1 \,\, {w_{lj}}^{-1} {E_{kj}}^{-1} {E_{ij}}^{-1} w_{lj}) \otimes e_k \\
   & \hspace{2em} + ({w_{lj}}^{-1} E_{ij} w_{lj} \,\, s_2 \,\, {w_{lj}}^{-1} {E_{ij}}^{-1} w_{lj}) \otimes e_k + s_3 \otimes e_k
\end{split}
\end{equation}
where
\[\begin{split}
   s_1 & := ({w_{lj}}^{-1} E_{ik^{-1}} w_{lj})^{-1} {E_{ik}}^{-1}, \hspace{1em}
   s_2 := ({w_{lj}}^{-1} E_{kj^{-1}} w_{lj})^{-1} {E_{kl}}^{-1}, \\
   s_3 & := ({w_{lj}}^{-1} E_{ij^{-1}} w_{lj})^{-1} {E_{il}}^{-1}. \\
  \end{split}\]
Then using Lemmas {\rmfamily \ref{L1}} and {\rmfamily \ref{L2}} repeatedly, we obtain
\begin{equation}
 r_{il}(k) \otimes e_k \equiv r_{ij}(k) \otimes e_k + s_1 \otimes e_l + s_3 \otimes e_k. \label{eq20}
\end{equation}
On the other hand, using Lemma {\rmfamily \ref{LB}}, we see
\[\begin{split}
    s_1 & = ({w_{lj}}^{-1} E_{ik^{-1}} w_{lj})^{-1} {E_{ik^{-1}}} \\
        & = E_{j^{-1}l} E_{l^{-1}j^{-1}} \, [E_{jl^{-1}},E_{ik}] \, E_{l^{-1}j} E_{j^{-1}l^{-1}} 
             \cdot E_{j^{-1}l} \, [E_{l^{-1}j^{-1}},E_{ik}] \, E_{j^{-1}l^{-1}}
        \cdot [E_{j^{-1}l},E_{ik}],
    \end{split}\]
and hence
\[ s_1 \otimes e_l \equiv r_{ik}(\cdot) \otimes e_l - r_{ik}(\cdot) \otimes e_j. \]
Furthermore, applying (vi) of Lemma {\rmfamily \ref{LA}} to $s_3$, we have
\[\begin{split}
    s_3 & = ({w_{lj}}^{-1} E_{ij^{-1}} w_{lj})^{-1} {E_{il^{-1}}} \\
        &= E_{j^{-1}l} E_{l^{-1}j^{-1}} E_{il} \, {r_{il^{-1}}(j)}^{-1} \, E_{il^{-1}} E_{l^{-1}j} E_{j^{-1}l^{-1}}
               \cdot E_{j^{-1}l} E_{il} \, {r_{ij^{-1}}(l^{-1})} \, E_{il^{-1}} E_{j^{-1}l^{-1}} \\
        & \hspace{2em} \cdot E_{j^{-1}l} E_{il} \, [E_{ij^{-1}},E_{l^{-1}i^{-1}}] \, E_{il^{-1}} E_{j^{-1}l^{-1}}
               \cdot [E_{j^{-1}l},E_{il}],
    \end{split}\]
and hence
\[ s_3 \otimes e_k \equiv r_{il}(\cdot) \otimes e_k - r_{ij}(\cdot) \otimes e_k. \]
Substituting these results into (\ref{eq20}), we obtain $S_{ijk} \equiv S_{ilk}$.

By the same argument, considering the equation
$(x_l^{-1},x_j,x_i,x_j,x_k) \otimes e_k$, we obtain
$S_{ijk} \equiv - S_{ilk}$, and $2 S_{ijk} \equiv 0$. Then $2$ is invertible in $L$,
we obtain $S_{ijk} \equiv 0$, i.e.,
\[ r_{ij}(k) \otimes e_k \equiv r_{ij}(\cdot) \otimes e_k + r_{ik}(\cdot) \otimes e_j. \]

Similarly, considering the equations $(x_k^{-1},x_l,x_i,x_j,x_l) \otimes e_k$ and
$(x_i^{-1},x_l,x_l,x_j,x_k^{\pm1}) \otimes e_k$,
we obtain
\[\begin{split}
   & r_{ij}(k^{-1}) \otimes e_k \equiv r_{ij}(\cdot) \otimes e_k + r_{ik}(\cdot) \otimes e_j, \\
   & r_{{i}^{-1}j}(k^{\pm1}) \otimes e_k \equiv r_{{i}^{-1}j}(\cdot) \otimes e_k + r_{{i}^{-1}k}(\cdot) \otimes e_j
  \end{split}\]
respectively. (For details, see \cite{Sa2}.)

\vspace{0.5em}

By the argument above, we can remove the generators $r_{{i}^{\pm1}j}(k^{\pm1}) \otimes e_k$
from the generationg set $\mathfrak{E}$.

\vspace{1em}

{\bf Step 4.} Here we consider the generators $h_{ij} \otimes e_p$.

\vspace{0.5em}

First we consider the case $p \neq i,j$.
From Lemma {\rmfamily \ref{LB}}, we have
\[\begin{split}
    [{w_{ij}}^{-1}, E_{lk}] & = ({w_{ij}}^{-1} E_{lk^{-1}} w_{ij})^{-1} {E_{lk^{-1}}} \\
        & = E_{j^{-1}i} E_{i^{-1}j^{-1}} \, [E_{ji^{-1}},E_{lk}] \, E_{i^{-1}j} E_{j^{-1}i^{-1}} \\
        & \hspace{2em}  \cdot E_{j^{-1}i} \, [E_{i^{-1}j^{-1}},E_{lk}] \, E_{j^{-1}i^{-1}}
           \cdot [E_{j^{-1}i},E_{lk}],
    \end{split}\]
and hence
\begin{equation}
[{w_{ij}}^{-1}, E_{lk}] \otimes e_l \equiv  -r_{ji}(\cdot) \otimes e_k 
- r_{i^{-1}j}(\cdot) \otimes e_k + r_{j^{-1}i}(\cdot) \otimes e_k. \label{eq3}
\end{equation}
On the other hand, observing (\ref{eq2}), we have
\[ [{w_{ij}}^{-1},  E_{lk}] = ({w_{ij}}^{-1} E_{lk} {h_{ij}}^{-1} E_{lk^{-1}} w_{ij} ) \cdot
 ({w_{ij}}^{-1} {h_{ij}} w_{ij} ) \cdot [w_{i^{-1}j}, E_{lk}].\] 
Using Lemmas {\rmfamily \ref{L1}} and {\rmfamily \ref{L2}}, we have
\[\begin{split}
   & ({w_{ij}}^{-1} {h_{ij}} w_{ij} ) \otimes e_l \equiv h_{ij} \otimes e_l, \\
   & ({w_{ij}}^{-1} E_{lk} {h_{ij}}^{-1} E_{lk^{-1}} w_{ij} ) \otimes e_l \equiv - h_{ij} \otimes e_l - h_{ij} \otimes e_k
  \end{split}\]
Furthermore, computing $[w_{i^{-1}j}, E_{lk}] \otimes e_l$ in a way similar to (\ref{eq3}), we have
\[ [w_{i^{-1}j}, E_{lk}] \otimes e_l \equiv r_{j^{-1}i}(\cdot) \otimes e_k + r_{ij}(\cdot) \otimes e_k - r_{ji}(\cdot) \otimes e_k. \]
Hence
\begin{equation}
[{w_{ij}}^{-1},  E_{lk}] \otimes e_l \equiv - h_{ij} \otimes e_k + r_{j^{-1}i}(\cdot) \otimes e_k + 
   r_{ij}(\cdot) \otimes e_k - r_{ji}(\cdot) \otimes e_k. \label{eq4}
\end{equation}
Comparing (\ref{eq3}) with (\ref{eq4}), we obtain
\[ h_{ij} \otimes e_k \equiv -r_{ij}(\cdot) \otimes e_k - r_{i^{-1}j}(\cdot) \otimes e_k. \]

\vspace{0.5em}

Next we consider the case $p=i$. 
Applying (iv) of Lemma {\rmfamily \ref{LA}} to
$\{ ({w_{ij}}^{-1} E_{ki^{-1}} w_{ij})^{-1} E_{kj} \} \otimes e_k$, we see
\[\begin{split}
    ({w_{ij}}^{-1} E_{ki^{-1}} w_{ij})^{-1} E_{kj}
        & = E_{j^{-1}i} E_{i^{-1}j^{-1}} \, [E_{ji^{-1}},E_{ki}] \, E_{i^{-1}j} E_{j^{-1}i^{-1}} \\
        & \hspace{2em} \cdot E_{j^{-1}i} E_{ki} \, {r_{kj^{-1}}(i^{-1})}^{-1} \, E_{ki^{-1}} E_{j^{-1}i^{-1}} \\
        & \hspace{2em} \cdot E_{j^{-1}i} E_{ki} \, {r_{ki^{-1}}(j^{-1})}^{-1} \, E_{ki^{-1}} E_{j^{-1}i^{-1}},
  \end{split}\]
and hence
\begin{equation}
\begin{split}
 \{ ({w_{ij}}^{-1} E_{ki^{-1}} w_{ij})^{-1} E_{kj} \} \otimes e_k & \equiv - r_{ji}(\cdot) \otimes e_i + r_{kj}(i^{-1}) \otimes e_k
           + r_{kj}(i^{-1}) \otimes e_i \\
   & \hspace{2em} + r_{ki}(j^{-1}) \otimes e_k + r_{ki}(\cdot) \otimes e_i.
\end{split} \label{heq1}
\end{equation}
On the other hand, using (\ref{eq2}), we have
\[\begin{split}
    ({w_{ij}}^{-1} E_{ki^{-1}} w_{ij})^{-1} E_{kj}
       & = ({w_{ij}}^{-1} E_{ki} h_{ij} E_{ki^{-1}} w_{ij}) \cdot ({w_{ij}}^{-1} {h_{ij}}^{-1} w_{ij}) \\
       & \hspace{2em} \cdot ({w_{i^{-1}j^{-1}}}^{-1} E_{ki^{-1}} w_{i^{-1}j^{-1}})^{-1} {E_{kj}}.
  \end{split} \]
Tensoring both hands side of the equation above with $e_k$, we have
\begin{equation}
 h_{ij} \otimes e_i \equiv \{({w_{ij}}^{-1} E_{ki^{-1}} w_{ij})^{-1} E_{kj}\} \otimes e_k
   - \{ ({w_{i^{-1}j^{-1}}}^{-1} E_{ki^{-1}} w_{i^{-1}j^{-1}})^{-1} {E_{kj}} \} \otimes e_k. \label{heq2}
\end{equation}
Applying (iii) of Lemma {\rmfamily \ref{LA}} to
$({w_{i^{-1}j^{-1}}}^{-1} E_{ki^{-1}} w_{i^{-1}j^{-1}})^{-1} {E_{kj}}$, we see
\[\begin{split}
   ({w_{i^{-1}j^{-1}}}^{-1} E_{ki^{-1}} w_{i^{-1}j^{-1}})^{-1} {E_{kj}}
        & = E_{ji^{-1}} E_{ij} \, [E_{j^{-1}i},E_{ki}] \, E_{ij^{-1}} E_{ji} \\
        & \hspace{2em} \cdot E_{ji^{-1}} E_{kj^{-1}} \, {r_{kj}(i)}^{-1} \,
                                  E_{kj} E_{ji} \\
        & \hspace{2em} \cdot E_{ji^{-1}} E_{kj^{-1}} \, {r_{ki}(j)}^{-1} \,
                                  E_{kj} E_{ji}
    \end{split} \]
and
\begin{equation}
\begin{split}
 \{ ({w_{i^{-1}j^{-1}}}^{-1} E_{ki^{-1}} w_{i^{-1}j^{-1}})^{-1} {E_{kj}} \} \otimes e_k & \equiv r_{j^{-1}i}(\cdot) \otimes e_i
     - r_{kj}(i) \otimes e_k + r_{kj}(\cdot) \otimes e_j \\
     & \hspace{2em} - r_{ki}(j) \otimes e_k + r_{ki}(j) \otimes e_j.
\end{split} \label{heq3}
\end{equation}
Substituiting (\ref{heq1}) and (\ref{heq3}) into (\ref{heq2}), we obtain
\[\begin{split}
   h_{ij} \otimes e_i & \equiv - r_{j^{-1}i}(\cdot) \otimes e_i - r_{ji}(\cdot) \otimes e_i
                               + r_{kj}(i) \otimes e_k + r_{kj}(i^{-1}) \otimes  e_k \\
                      & \hspace{1em} + r_{ki}(j) \otimes e_k + r_{ki}(j^{-1}) \otimes  e_k
                               - r_{kj}(\cdot) \otimes e_j + r_{ki}(\cdot) \otimes e_i \\
                      & \hspace{1em} - r_{ki}(j) \otimes e_j + r_{kj}(i^{-1}) \otimes e_i
  \end{split}\]
Similarly, considering
$\{ ({w_{ij}} E_{ki^{-1}} w_{ij^{-1}})^{-1}$ $E_{kj^{-1}} \} \otimes e_k$, we have
\[\begin{split}
   h_{ij} \otimes e_j & \equiv r_{j^{-1}i}(\cdot) \otimes e_i + r_{ji}(\cdot) \otimes e_i
                               + r_{kj}(i^{-1}) \otimes e_k + r_{kj}(i) \otimes  e_k \\
                      & \hspace{1em} - r_{ki}(j) \otimes e_k - r_{ki}(j^{-1}) \otimes  e_k
                               - r_{ki}(\cdot) \otimes e_i + r_{kj}(\cdot) \otimes e_j \\
                      & \hspace{1em} + r_{kj}(i^{-1}) \otimes e_j - r_{ki}(j^{-1}) \otimes e_j
  \end{split}\]

By the argument above, we can remove the generators $h_{ij} \otimes e_p$
from the generationg set $\mathfrak{E}$.

\vspace{1em}

{\bf Step 5.} Here we consider the generators $r_{{i}^{\pm1}j}(k^{\pm1}) \otimes e_i$.

\vspace{0.5em}

For convenience, we use the following notation.
Let $V$ be the quotient $L$-module of $\bar{R}^{\mathrm{ab}} {\otimes}_{\mathrm{Aut}^{+}F_n} H_L$ by
the $L$-submodule generated by the elements $r_{i^{\pm1}j}(\cdot) \otimes e_k$ for $k \neq i$.
We use $\c$ for the equality in $V$.

First we consider the equation $(x_l,x_k,x_i,x_j,x_k) \otimes e_i$ for distinct $i,j,k$ and $l$.
It is given by
\[\begin{split}
   r_{ij}(l) \otimes e_i & =
     s_1^{-1} \otimes e_i + ({w_{lk}}^{-1} E_{ik} w_{lk} \,\,
     s_2^{-1} \,\, {w_{lk}}^{-1} {E_{ik}}^{-1} w_{lk}) \otimes e_i \\
   & \hspace{2em} + ({w_{lk}}^{-1} r_{ij}(k) w_{lk}) \otimes e_i \\
   & \hspace{2em} + ({w_{lk}}^{-1} E_{ij} E_{kj} w_{lk} \,\, s_1
      \,\, {w_{lk}}^{-1} E_{kj}^{-1} {E_{ij}}^{-1} w_{lk}) \otimes e_i \\
   & \hspace{2em} + ({w_{lk}}^{-1} E_{ij} w_{lk} \,\, s_2
      \,\, {w_{lk}}^{-1} {E_{ij}}^{-1} w_{lk}) \otimes e_i + s_3 \otimes e_i
\end{split}\]
where
\[\begin{split}
   s_1 & := ({w_{lk}}^{-1} E_{ik^{-1}} w_{lk})^{-1} {E_{il}}^{-1}, \hspace{1em}
   s_2 := ({w_{lk}}^{-1} E_{kj^{-1}} w_{lk})^{-1} {E_{lj}}^{-1}, \\
   s_3 & := ({w_{lk}}^{-1} E_{ij^{-1}} w_{lk})^{-1} {E_{ij}}^{-1}. \\
  \end{split}\]
Then using Lemmas {\rmfamily \ref{L1}} and {\rmfamily \ref{L2}}, we obtain
\begin{equation}
  r_{ij}(l) \otimes e_i \equiv  r_{ij}(k) \otimes e_i + s_1 \otimes e_j + s_2 \otimes e_j
             - s_2 \otimes e_l + s_3 \otimes e_i. \label{eq60}
\end{equation}
By an argumet similar to that in {\bf{(3-b)}}, we can compute
\[ s_1 \otimes e_j \equiv r_{il}(\cdot) \otimes e_j - r_{ik}(\cdot) \otimes e_j \c 0 \]
and
\[ s_3 \otimes e_i \equiv -r_{kl}(\cdot) \otimes e_j - r_{l^{-1}k}(\cdot) \otimes e_j + r_{k^{-1}l}(\cdot) \otimes e_j
    \c 0. \]
On the other hand, using (\ref{eq2}), we have
\[\begin{split}
   s_2 = ({w_{lk}}^{-1} E_{kj} h_{lk} E_{kj^{-1}} w_{lk}) \cdot ({w_{lk}}^{-1} {h_{lk}}^{-1} w_{lk})
          \cdot ({w_{l^{-1}k^{-1}}}^{-1} E_{kj^{-1}} w_{l^{-1}k^{-1}})^{-1} {E_{lj^{-1}}},
  \end{split}\]
and hence
\[\begin{split}
   s_2 \otimes e_j & \equiv \{ ({w_{l^{-1}k^{-1}}}^{-1} E_{kj^{-1}} w_{l^{-1}k^{-1}})^{-1} {E_{lj^{-1}}} \} \otimes e_j, \\
   s_2 \otimes e_l & \equiv \{ ({w_{l^{-1}k^{-1}}}^{-1} E_{kj^{-1}} w_{l^{-1}k^{-1}})^{-1} {E_{lj^{-1}}} \} \otimes e_l
                      + h_{lk} \otimes e_j, \\
                   & \equiv \{ ({w_{l^{-1}k^{-1}}}^{-1} E_{kj^{-1}} w_{l^{-1}k^{-1}})^{-1} {E_{lj^{-1}}} \} \otimes e_l
                      - r_{lk}(\cdot) \otimes e_j - r_{l^{-1}k}(\cdot) \otimes e_j.
  \end{split}\]
Then, applying (ii) of Lemma {\rmfamily \ref{LA}} to
$({w_{l^{-1}k^{-1}}}^{-1} E_{kj^{-1}} w_{l^{-1}k^{-1}})^{-1} {E_{lj^{-1}}}$, we see
\[\begin{split}
    ({w_{l^{-1}k^{-1}}}^{-1} E_{kj^{-1}} w_{l^{-1}k^{-1}})^{-1} {E_{lj^{-1}}}
        & = E_{kl^{-1}} E_{lk} \, [E_{k^{-1}l},E_{kj}] \, E_{lk^{-1}} E_{kl} \\
        & \hspace{2em} \cdot E_{kl^{-1}} \, {r_{lj}(k)} \, E_{kl} \cdot E_{lj} E_{kl^{-1}} \, r_{kj^{-1}}(l) \, E_{kl} E_{lj^{-1}},
   \end{split}\]
and hence
\[\begin{split}
   s_2 \otimes e_j & \equiv r_{lj}(\cdot) \otimes e_j - r_{kj}(\cdot) \otimes e_j \c 0, \\
   s_2 \otimes e_l & \equiv r_{kj}(\cdot) \otimes e_l + r_{k^{-1}l}(\cdot) \otimes e_j + r_{lj}(k) \otimes e_l \\
                   & \hspace{4em}  - r_{kj}(l) \otimes e_l - r_{kj}(\cdot) \otimes e_j
                                   - r_{lk}(\cdot) \otimes e_j - r_{l^{-1}k}(\cdot) \otimes e_j, \\
                   & \equiv r_{kj}(\cdot) \otimes e_l + r_{k^{-1}l}(\cdot) \otimes e_j + r_{lj}(k) \otimes e_l \\
                   & \hspace{4em}  - r_{kj}(\cdot) \otimes e_l - r_{kl}(\cdot) \otimes e_j - r_{kj}(\cdot) \otimes e_j
                                   - r_{lk}(\cdot) \otimes e_j - r_{l^{-1}k}(\cdot) \otimes e_j, \\
                   & \c r_{lj}(k) \otimes e_l.
  \end{split}\]
Substituting these results into (\ref{eq60}), we obtain
\begin{equation}
r_{ij}(l) \otimes e_i \c r_{ij}(k) \otimes e_i - r_{lj}(k) \otimes e_l. \label{eq6}
\end{equation}

Similarly, considering the equations $(x_l^{-1},x_i,x_i,x_j,x_k) \otimes e_l$, $(x_k,x_l,x_i,x_j,x_k^{-1}) \otimes e_i$
and $(x_l^{-1},x_k,x_i^{-1},x_j,x_k^{-1}) \otimes e_i$, we obtain
\begin{gather}
r_{ij}(k) \otimes e_i \c r_{l^{-1}j}(i) \otimes e_l - r_{l^{-1}j}(k) \otimes e_l, \label{eq7} \\ 
r_{ij}(k^{-1}) \otimes e_i \c r_{ij}(l) \otimes e_i - r_{k^{-1}j}(l) \otimes e_k, \notag \\ 
r_{i^{-1}j}(k^{-1}) \otimes e_i \c r_{i^{-1}j}(l) \otimes e_i - r_{lj}(k^{-1}) \otimes e_l \notag
\end{gather}

From the equations above, we see that $V$ is generated by $r_{i^{\pm1}j}(k) \otimes e_i$.
We reduce these generators of $V$ more.
On the equation (\ref{eq6}), exchanging the roles of $k$ and $l$,
we obtain
\[ r_{ij}(k) \otimes e_i \c r_{ij}(l) \otimes e_i - r_{kj}(l) \otimes e_k, \]
and hence
\[ r_{kj}(l) \otimes e_k \c - r_{lj}(k) \otimes e_l. \]
For any $j \in \{ 1, \ldots , n \}$,
choose a number $\mu_j \in \{ 1, \ldots , n \}$ such that $\mu_j \neq j$ and fix it.
Then we have
\[ r_{\mu_j j}(k) \otimes e_{\mu_j} \c - r_{kj}(\mu_j) \otimes e_k,
     \hspace{1em} r_{ij}(k) \otimes e_i \c r_{ij}(\mu_j) \otimes e_i - r_{kj}(\mu_j) \otimes e_k. \]
Furthermore, from (\ref{eq7}), we have
\[ r_{i^{-1}j}(k) \otimes e_i \c r_{kj}(\mu_j) \otimes e_k + r_{i^{-1}j}(\mu_j) \otimes e_i. \]
This shows that the $L$-module $V$ is generated by
\[ r_{\alpha j}(\mu_j) \otimes e_{\alpha}, \hspace{1em} (1 \leq \alpha \leq n, \,\,\, \alpha \neq j, \mu_j) \]
and
\[ r_{{\beta}^{-1} j}(\mu_j) \otimes e_{\beta}, \hspace{1em} (1 \leq \beta \leq n, \,\,\, \beta \neq j). \]

\vspace{1em}

Therefore we conclude that the generating set $\mathfrak{E}$ of $\bar{R}^{\mathrm{ab}} {\otimes}_{\mathrm{Aut}^{+}F_n} H_L$
is reduced to
\[ \{ r_{i^{\pm1}j}(\cdot) \otimes e_p \,|\, p \neq i \} \, \cup \, 
   \{ r_{\alpha j}(\mu_j) \otimes e_{\alpha} \,|\, 1 \leq j \leq n \} \,
   \cup \, \{ r_{\beta^{-1} j}(\mu_j) \otimes_{\beta} \,|\, 1 \leq \beta \leq n \}. \]
The number of the generators above is just $2n(n^2-n)-n$.
This completes the proof of Proposition {\rmfamily \ref{P1}}. $\square$

\subsection{The proof of Proposotion {\rmfamily \ref{P2}}}\label{H^*}
\hspace*{\fill}\ 

\vspace{0.5em}

In this subsection, we prove Proposotion {\rmfamily \ref{P2}}.
The outline of the proof is similar to that of Proposotion {\rmfamily \ref{P1}}.
Since the image of the map $H_1(\bar{R},H_L^*)_{\mathrm{Aut}^{+}F_n} \rightarrow H_1(F,H_L^*)$ is
isomorphic to the free $L$-module of rank $2n(n^2-n)-n-1$,
$H_1(\bar{R},H_L^*)_{\mathrm{Aut}^{+}F_n}$ contains a free $L$-submodule which rank is
greater than or equal to $2n(n^2-n)-n-1$.
To show it is just $2n(n^2-n)-n-1$, it suffices to show that
$H_1(\bar{R},H_L)_{\mathrm{Aut}^{+}F_n}$ is generated by just $2n(n^2-n)-n-1$ elements.
We have $H_1(\bar{R},H_L^*)_{\mathrm{Aut}^{+}F_n} = \bar{R}^{\mathrm{ab}} {\otimes}_{\mathrm{Aut}^{+}F_n} H_L^*$,
and see that
\[ \mathfrak{E}^* := \{ r \otimes e_p^* \,|\, r \in R, \,\, 1 \leq p \leq n \} \]
is a generating set of $\bar{R}^{\mathrm{ab}} {\otimes}_{\mathrm{Aut}^{+}F_n} H_L^*$. In the following, we reduce the elemets
of $\mathfrak{E}^*$.
We also use $\equiv$ for the equality in $\bar{R}^{\mathrm{ab}} {\otimes}_{\mathrm{Aut}^{+}F_n} H_L^*$.

\vspace{1em}

{\bf Step 0.} By an argument similar to that of Step 0 in Subsection {\rmfamily \ref{H}}, we have
\begin{lem}\label{L1*}
For $n \geq 3$,
\[\begin{split}
   (E_{{i}^{\pm 1}j} \, r \, {E_{{i}^{\pm 1}j^{-1}}}) \otimes e_p^* & \equiv
              \begin{cases}
                r \otimes e_p^*, & \hspace{1em} p \neq j, \\
                r \otimes e_j^* \mp r \otimes e_i^*, & \hspace{1em} p=j,
              \end{cases} \\
   ({E_{{i}^{\pm 1}j^{-1}}} \, r \, E_{{i}^{\pm 1}j}) \otimes e_p^* & \equiv
              \begin{cases}
                r \otimes e_p^*, & \hspace{1em} p \neq j, \\
                r \otimes e_j^* \pm r \otimes e_i^*, & \hspace{1em} p=j.
              \end{cases} \\
  \end{split}\]
\end{lem}

\begin{cor}\label{C1*}
For $n \geq 3$,
\[ [E_{{i}^{\pm1} j}, r] \otimes e_p^* \equiv
              \begin{cases}
                0, & \hspace{1em} p \neq j, \\
                \mp r \otimes e_i^*, & \hspace{1em} p=j,
              \end{cases}\]
\[ [{E_{{i}^{\pm1} j^{-1}}}, r] \otimes e_p^* \equiv 
              \begin{cases}
                0, & \hspace{1em} p \neq j, \\
                \pm r \otimes e_i^*, & \hspace{1em} p=j.
              \end{cases}\]
\end{cor}

\begin{lem}\label{L2*}
For $n \geq 3$,
\[(w_{{i}^{\pm 1}j} \, r \, {w_{{i}^{\pm 1}j}}^{-1}) \otimes e_p^* \equiv
              \begin{cases}
                r \otimes e_p^*, & \hspace{1em} p \neq i,j, \\
                \mp r \otimes e_j^*, & \hspace{1em} p=i, \\
                \pm r \otimes e_i^*, & \hspace{1em} p=j,
              \end{cases}\]
\[({w_{{i}^{\pm 1}j}}^{-1} \, r \, {w_{{i}^{\pm 1}j}}) \otimes e_p^* \equiv
              \begin{cases}
                r \otimes e_p^*, & \hspace{1em} p \neq i,j, \\
                \pm r \otimes e_j^*, & \hspace{1em} p=i, \\
                \mp r \otimes e_i^*, & \hspace{1em} p=j.
              \end{cases}\]
\end{lem}

\begin{cor}\label{C2*}
For $n \geq 3$,
\[ [w_{{i}^{\pm1}j}, r] \otimes e_p^* \equiv
   \begin{cases}
     0,  \hspace{1em} & p \neq i,j, \\
     -r \otimes e_i^* \mp r \otimes e_j^*, \hspace{1em} & p = i, \\
     \pm r \otimes e_i^* -r \otimes e_j^*, \hspace{1em} & p = j, \\
   \end{cases}\]
\[ [w_{{i}^{\pm1}j^{-1}}, r] \otimes e_p^* \equiv
   \begin{cases}
     0,  \hspace{1em} & p \neq i,j, \\
     -r \otimes e_i^* \pm r \otimes e_j^*, \hspace{1em} & p = i, \\
     \mp r \otimes e_i^* -r \otimes e_j^*, \hspace{1em} & p = j. \\
   \end{cases}\]
\end{cor}

\vspace{0.5em}

Considering any relator of (R2) of the Gersten's presentation is conjugate to one of the relator of (R2-1), $\ldots$ , (R2-8),
or considering Lemma {\rmfamily \ref{LC}}, for any relator $r=$ (R2), (R3) and (R4),
we can rewrite a element $r \otimes e_p^*$ with the relators (R2-1), $\ldots$ , (R4-1)
using Lemmas {\rmfamily \ref{L1}} and {\rmfamily \ref{L2}}.

\vspace{1em}

{\bf Step 1.} First we consider the generators ${w_{ij}}^4 \otimes e_p^*$.
By the same argument as that of Step 1 in Subsection {\rmfamily \ref{H}}, we see
\[ {w_{ij}}^4 \otimes e_p^* = \frac{1}{2} \, 2 \, {w_{ij}}^4 \otimes e_p^* \equiv \frac{1}{2} \, {w_{ij}}^8 \otimes e_p^* \]
is rewritten as a sum of the generators $r \otimes e_p^*$ for $r=$ (R2-1), $\ldots$ , (R4-1).
Therefore we can remove the generators ${w_{ij}}^4 \otimes e_p^*$ from the generating set $\mathfrak{E}^*$.

\vspace{1em}

{\bf Step 2.}
Here we show that the generators $r \otimes e_p^*$ for $r=$ (R2-1), $\ldots$ , (R2-8)
is zero or equal to one of the generators $r_{{i}^{\pm1}j}(k^{\pm1}) \otimes e_p^*$. We have
\begin{lem}\label{L3*}
For $n \geq 6$ and distinct $i,j,k,l$ and $m$, we have
\begin{enumerate}
\item[(i)] {\bf{(R2-6):}}
\[ [E_{ij}, E_{kl}] \otimes e_p^* \equiv
   \begin{cases}
     0,  \hspace{1em} & p \neq j,l, \\
     r_{kl}(m) \otimes e_i^*, \hspace{1em} & p = j, \\
     - r_{ij}(m) \otimes e_k^*, \hspace{1em} & p = l. \\
   \end{cases}\]
\item[(ii)] {\bf{(R2-7), (R2-8):}}
\[ [E_{{i}^{-1}j}, E_{{k}^{\pm1}l}] \otimes e_p^* \equiv
   \begin{cases}
     0,  \hspace{1em} & p \neq j,l \\
     - r_{{k}^{\pm1}l}(m) \otimes e_i^*, \hspace{1em} & p = j, \\
     \mp r_{{i}^{-1}j}(m) \otimes e_k^*, \hspace{1em} & p = l. \\
   \end{cases}\]
\item[(iii)] {\bf{(R2-5):}}
\[ [E_{ij}, E_{{i}^{-1}k}] \otimes e_p^* \equiv
   \begin{cases}
     0,  \hspace{1em} & p \neq j,k, \\
     r_{{i}^{-1}k}(l) \otimes e_i^*, \hspace{1em} & p = j, \\
     r_{ij}(l) \otimes e_i, \hspace{1em} & p = k.
   \end{cases}\]
\item[(iv)] {\bf{(R2-1):}}
\[ [E_{ij}, E_{{i}^{-1}j}] \otimes e_p^* \equiv
   \begin{cases}
     0,  \hspace{1em} & p \neq j, \\
     r_{ij}(k) \otimes e_i^* - r_{{i}^{-1}j}(l) \otimes e_i, \hspace{1em} & p = j.
   \end{cases}\]
\item[(v)] {\bf{(R2-2):}}
\[ [E_{ij}, E_{kj}] \otimes e_p^* \equiv
   \begin{cases}
     0,  \hspace{1em} & p \neq j, \\
     r_{kj}(l) \otimes e_i^* -r_{ij}(m) \otimes e_k^*, \hspace{1em} & p = j.
   \end{cases}\]
\item[(vi)] {\bf{(R2-3), (R2-4):}}
\[ [E_{{i}^{-1}j}, E_{{k}^{\pm1}j}] \otimes e_p^* \equiv
   \begin{cases}
     0,  \hspace{1em} & p \neq j, \\
     - r_{{k}^{\pm1}j}(l) \otimes e_i^* \mp r_{{i}^{-1}j}(m) \otimes e_k^*, \hspace{1em} & p = j. \\
   \end{cases}\]
\end{enumerate}
\end{lem}
Since this Lemma is proved by an argument similar to that in Lemma {\rmfamily \ref{L3}}, we omit the details.
(For details, see \cite{Sa2}.)

\vspace{0.5em}

By the lemma above, we can remove the generators $r \otimes e_p^*$ for $r=$ (R2-1), $\ldots$ , (R2-8)
from the generationg set $\mathfrak{E}^*$.

\vspace{1em}

{\bf Step 3.} Here we consider the generators $r_{{i}^{\pm1}j}(k^{\pm1}) \otimes e_p^*$ for $p \neq j$.

\vspace{0.5em}

{\bf (3-a)} {\bf The case $p \neq j,k$}.

\vspace{0.5em}

First we consider the case $p \neq k$. By an argument similar to that of {\bf (3-a)} in Subsection {\rmfamily \ref{H}},
observing the results of Lemma {\rmfamily \ref{L3*}}, we can set
\[\begin{split}
   r_{ij}(\cdot) \otimes e_p^* :\equiv & r_{ij}(k^{\pm1}) \otimes e_p^*, \\
   r_{{i}^{-1}j}(\cdot) \otimes e_p^* :\equiv & r_{{i}^{-1}j}(k^{\pm1}) \otimes e_p^*
  \end{split}\]
for $p \neq j,k$.

\vspace{1em}

{\bf (3-b)} {\bf The case $p=k$.}

\vspace{0.5em}

For the case $p=k$, set
\[ S_{ijk}^* := r_{ij}(k) \otimes e_k^* - r_{ij}(\cdot) \otimes e_k^* - r_{kj}(\cdot) \otimes e_i^*. \]
By the same argument as that of {\bf (3-b)} in Subsection {\rmfamily \ref{H}},
observing the equations $(x_l^{\pm1},x_j, x_i,x_j,x_k) \otimes e_k^*$,
we obtain $S_{ijk}^* \equiv 0$.
Furthermore, observing $(x_k^{-1},x_l,x_i,x_j,x_l) \otimes e_k^*$ and
$(x_i^{-1},x_l,x_l,x_j,x_k^{\pm1}) \otimes e_k^*$, we obtain
\[\begin{split}
   & r_{ij}(k^{-1}) \otimes e_k^* \equiv r_{ij}(\cdot) \otimes e_k^* - r_{k^{-1}j}(\cdot) \otimes e_i^*, \\
   & r_{{i}^{-1}j}(k^{\pm1}) \otimes e_k^* \equiv r_{{i}^{-1}j}(\cdot) \otimes e_k^* \mp r_{{k}^{\pm1}j}(\cdot) \otimes e_i^*.
  \end{split}\]

\vspace{0.5em}

By the argument above, we can remove the generators $r_{{i}^{\pm1}j}(k^{\pm1}) \otimes e_k^*$
from the generationg set $\mathfrak{E}^*$.

\vspace{1em}

{\bf Step 4.}
Here we consider the generators $h_{ij} \otimes e_p^*$ for $p \neq i,j$.
By an argument similar to that of Step 4 in Subsection {\rmfamily \ref{H}},
considering the elements $[{w_{ij}}^{-1},  E_{kl}] \otimes e_l^*$, we obtain
\[ h_{ij} \otimes e_k^* \equiv r_{ij}(\cdot) \otimes e_k^* + r_{i^{-1}j}(\cdot) \otimes e_k^*. \]
The cases where $p=i$ or $j$ are mensioned in Step 6 later.

\vspace{1em}

{\bf Step 5.}
Let $V$ be the quotient $L$-module of $\bar{R}^{\mathrm{ab}} {\otimes}_{\mathrm{Aut}^{+}F_n} H_L^*$ by
the $L$-submodule generated by the elements $r_{i^{\pm1}j}(\cdot) \otimes e_k^*$ for $k \neq j$.
Then from the argument above, the elements $r_{i^{\pm1}j}(k^{\pm1}) \otimes e_j^*$, $h_{ij} \otimes e_i^*$ and
$h_{ij} \otimes e_j^*$ generate $V$.
Here we reduce these generators of $V$. We use $\c$ for the equality in $V$.

First, considering the equation $(x_l,x_k,x_i,x_j,x_k) \otimes e_j^*$ in a way similar to that of Step 5
in Subsection {\rmfamily \ref{H}}, we have
\begin{equation}
r_{ik}(l^{-1}) \otimes e_k^* \c - r_{ij}(k) \otimes e_j^* + r_{ij}(l) \otimes e_j^*. \label{eq9}
\end{equation}
Similarly, considering $(x_l^{-1},x_j,x_i,x_k,x_j) \otimes e_k^*$, $(x_l,x_k,x_i^{-1},x_j,x_k) \otimes e_j^*$ and
$(x_l^{-1},x_j,x_i^{-1},x_k,x_j) \otimes e_k^*$, we obtain
\begin{gather}
r_{ik}(l^{-1}) \otimes e_k^* \c  r_{ik}(j) \otimes e_k^* + r_{ij}(l) \otimes e_j^*, \label{eq10} \\
r_{i^{-1}k}(l^{-1}) \otimes e_k^* \c - r_{i^{-1}j}(k) \otimes e_j^* + r_{i^{-1}j}(l) \otimes e_j^*, \label{eq22} \\
r_{i^{-1}k}(l^{-1}) \otimes e_k^* \c r_{i^{-1}k}(j) \otimes e_k^* + r_{i^{-1}j}(l) \otimes e_j^* \label{eq23}
\end{gather}
respectively.
Hence we see that the $L$-module $V$ is generated by $r_{i^{\pm1}j}(k) \otimes e_j^*$, $h_{ij} \otimes e_i^*$ and
$h_{ij} \otimes e_j^*$.

Substituting (\ref{eq10}) into (\ref{eq9}), and substituting (\ref{eq23}) into (\ref{eq22}), we obtain
\begin{equation}
r_{i^{\pm1}k}(j) \otimes e_k^* \c - r_{i^{\pm1}j}(k) \otimes e_j^*. \label{eq24}
\end{equation}
On the other hand, considering the equation $(x_l,x_k,x_i,x_j,x_k^{-1}) \otimes e_j^*$, we obtain
\[ r_{ij}(l^{-1}) \otimes e_j^* \c r_{ij}(k^{-1}) \otimes e_j^* + r_{ik}(l^{-1}) \otimes e_k^*. \]
Hence, rewriting each term of the equation above as a sum of
$r_{\alpha \beta}({\gamma}) \otimes e_{\beta}^*$, ($1 \leq \alpha, \beta, \gamma \leq n$), using (\ref{eq9}),
and using (\ref{eq24}), we have
$2(r_{ij}(k) \otimes e_j^* + r_{ik}(l) \otimes e_k^* - r_{ij}(l) \otimes e_j^*) \c 0$. Since $2$ is invertible in $V$, we obtain
\begin{equation}
 r_{ij}(k) \otimes e_j^* + r_{ik}(l) \otimes e_k^* - r_{ij}(l) \otimes e_j^* \c 0. \label{eq100}
\end{equation}
Similarly, considering $(x_i,x_l^{-1},x_l,x_j,x_k) \otimes e_j^*$, we have
\[\begin{split}
     r_{i^{-1}j}(k) \otimes e_j^* & \c r_{lj}(k) \otimes e_j^* - r_{i^{-1}k}(l) \otimes e_k^* \\
     & \hspace{2em} + r_{lk}(i^{-1}) \otimes e_k^* + r_{i^{-1}j}(l) \otimes e_j^* - r_{lj}(i^{-1}) \otimes e_j^*. 
  \end{split}\]
Using (\ref{eq9}) and (\ref{eq100}), we can reduce the equation above to
\[ r_{i^{-1}j}(k) \otimes e_j^* + r_{i^{-1}k}(l) \otimes e_k^* + r_{i^{-1}l}(j) \otimes e_l^* \c 0. \]

Now, using the equations above, we show that each generator $r_{i^{\pm1}j}(k) \otimes e_j^*$ is rewritten as a sum of
the generators type of $r_{i^{\pm1}j}(1) \otimes e_j^*$ and $r_{1^{\pm1}j}(2) \otimes e_j^*$.
For distinct $i, j, k \neq 1$, we have $r_{ij}(k) \otimes e_j^* \c - r_{ik}(1) \otimes e_k^* + r_{ij}(1) \otimes e_j^*$.
If $j=1$, we have $r_{i1}(k) \otimes e_1^* \c - r_{ik}(1) \otimes e_k^*$.
If $i =1$ and $j,k \neq 2$, then $r_{1j}(k) \otimes e_j^* \c - r_{1k}(2) \otimes e_k^* + r_{ij}(2) \otimes e_j^*$.
Finally, if $i=1$ and $j=2$, we have $r_{12}(k) \otimes e_2^* \c - r_{1k}(2) \otimes e_k^*$.
Hence any generator $r_{ij}(k) \otimes e_j^*$ is rewritten as a sum of the generators
$r_{ij}(1) \otimes e_j^*$ and $r_{1j}(2) \otimes e_j^*$.
Similarly we see that $r_{i^{-1}j}(k) \otimes e_j^*$ is rewritten as a sum of the generators
$r_{i^{-1}j}(1) \otimes e_j^*$ and $r_{1^{-1}j}(2) \otimes e_j^*$.

\vspace{0.5em}

From the argument above, we see that $V$ is generated by $r_{i^{\pm1}j}(1) \otimes e_j^*$, $r_{1^{\pm1}j}(2) \otimes e_j^*$,
$h_{ij} \otimes e_i^*$ and $h_{ij} \otimes e_j^*$, and hence
$\bar{R}^{\mathrm{ab}} {\otimes}_{\mathrm{Aut}^{+}F_n} H_L^*$ is generated by these elements and
$r_{i^{\pm1}j}(\cdot) \otimes e_k^*$.

\vspace{1em}

{\bf Step 6.}
Finally we consider the generators $h_{ij} \otimes e_p^*$ for $p = i,j$. 
Let $V'$ be the quotient $L$-module of $\bar{R}^{\mathrm{ab}} {\otimes}_{\mathrm{Aut}^{+}F_n} H_L^*$ by
the $L$-submodule generated by the elements $r_{i^{\pm1}j}(\cdot) \otimes e_k^*$, $r_{i^{\pm1}j}(1) \otimes e_j^*$ and
$r_{1^{\pm1}j}(2) \otimes e_j^*$. We use $\d$ for the equality in $V'$.

For distinct $i,j$ and $k$, the equation $\{ x_k,x_i^{-1},x_i,x_j \} \otimes e_j^*$ is given by
\[\begin{split}
   h_{k^{-1}j} \otimes e_j^* & =
     t_3 \otimes e_j^* + ({w_{ki^{-1}}}^{-1} E_{ji} w_{ki^{-1}} \,\, t_2
                            \,\, {w_{ki^{-1}}}^{-1} {E_{ji}}^{-1} w_{ki^{-1}}) \otimes e_j^* \\
   & \hspace{2em} + ({w_{ki^{-1}}}^{-1} E_{ji} E_{i^{-1}j} w_{ki^{-1}} \,\, t_1
                            \,\, {w_{ki^{-1}}}^{-1} {E_{i^{-1}j}}^{-1} {E_{ji}}^{-1} w_{ki^{-1}}) \otimes e_j^* \\
   & \hspace{2em} + ({w_{ki^{-1}}}^{-1} h_{ij} w_{ki^{-1}}) \otimes e_j^* \\
   & \hspace{2em} + ({w_{ki^{-1}}}^{-1} {E_{j^{-1}i}}^{-1} {E_{ij}}^{-1} w_{ki^{-1}}
                             \,\, t_4 \,\, {w_{ki^{-1}}}^{-1} E_{ij} {E_{j^{-1}i}} w_{ki^{-1}}) \otimes e_j^* \\
   & \hspace{2em} + ({w_{ki^{-1}}}^{-1} {E_{j^{-1}i}}^{-1} w_{ki^{-1}}
                             \,\, t_5 \,\, {w_{ki^{-1}}}^{-1} {E_{j^{-1}i}} w_{ki^{-1}}) \otimes e_j^* + t_6 \otimes e_j^*
  \end{split}\]
where
\[\begin{split}
   t_1 & :=  E_{j^{-1}k} (w_{ki^{-1}}^{-1} E_{j^{-1}i^{-1}} w_{ki^{-1}})^{-1},
                  \hspace{2.8em} t_2 := E_{k j} (w_{ki^{-1}}^{-1} E_{i^{-1}j} w_{ki^{-1}})^{-1}, \\
   t_3 & := E_{jk^{-1}} (w_{ki^{-1}}^{-1} E_{ji} w_{ki^{-1}})^{-1},
                        \hspace{7.25em} t_4 := (w_{ki^{-1}}^{-1} E_{ji^{-1}} w_{ki^{-1}})^{-1} E_{jk} , \\
   t_5 & := (w_{ki^{-1}}^{-1} E_{ij} w_{ki^{-1}})^{-1} E_{k^{-1} j},
                        \hspace{7.25em} t_6 := (w_{ki^{-1}}^{-1} E_{j^{-1}i} w_{ki^{-1}})^{-1} E_{j^{-1}k^{-1}}.
  \end{split}\]
Observing Lemma {\rmfamily \ref{LA}}, we see that
all $t_m$, ($1 \leq m \leq 6$), except for $t_2$ belong to the normal closure of the relators of (R2-1), $\ldots$ , (R3-4).
Hence, using Lemmas {\rmfamily \ref{L1*}} and Lemmas {\rmfamily \ref{L2*}}, we obtain
\begin{equation}
h_{k^{-1}j} \otimes e_j^* \d t_2 \otimes e_j^* + h_{ij} \otimes e_j^*. \label{eq50}
\end{equation}
From (\ref{eq2}), we have
\[\begin{split}
    t_2^{-1} & = ({w_{ki^{-1}}}^{-1} E_{i^{-1}j^{-1}}  w_{ki^{-1}})^{-1} {E_{kj^{-1}}} \\
             & = ({w_{ki^{-1}}}^{-1} E_{i^{-1}j} h_{ki^{-1}} E_{i^{-1}j^{-1}} w_{ki^{-1}}) 
                     \cdot ({w_{ki^{-1}}}^{-1} {h_{ki^{-1}}}^{-1} w_{ki^{-1}}) \\
             & \hspace{2em} \cdot ({w_{k^{-1}i}}^{-1} E_{i^{-1}j^{-1}} w_{k^{-1}i})^{-1} {E_{kj^{-1}}},
  \end{split} \]
and hence
\[  t_2^{-1} \otimes e_j^* \equiv \{ ({w_{k^{-1}i}}^{-1} E_{i^{-1}j^{-1}} w_{k^{-1}i})^{-1} {E_{kj^{-1}}} \} \otimes e_j^*
            - h_{ki^{-1}} \otimes e_i^*. \]
On the other hand, from Lemma {\rmfamily \ref{LC}}, we have
\[ h_{k^{-1}j} \otimes e_j^* \equiv ({w_{kj}}^{-1} {h_{kj}}^{-1} w_{kj}) \otimes e_j^* \equiv - h_{kj} \otimes e_k^*, \]
\[ h_{ki^{-1}} \otimes e_i^* \equiv ({w_{ki}}^{-1} {h_{ki}}^{-1} w_{ki}) \otimes e_i^* \equiv - h_{ki} \otimes e_k^*, \]
and see that the element $\{ ({w_{k^{-1}i}}^{-1} E_{i^{-1}j^{-1}} w_{k^{-1}i})^{-1} {E_{kj^{-1}}} \} \in \bar{R}$
belongs to the normal closure of the relators of (R2-1), $\ldots$ , (R3-4) by Lemma {\rmfamily \ref{LA}}.
Therefore we obtain
\[ t_2 \otimes e_j^* \d h_{ki} \otimes e_k^*. \]
Substituting these results into (\ref{eq50}), we obtain
\begin{equation}
h_{ij} \otimes e_j^* \d -h_{ki} \otimes e_k^* - h_{kj} \otimes e_k^*. \label{eq11}
\end{equation}
Similarly, considering $\{ x_k,x_i,x_i,x_j \} \otimes e_j^*$ and $\{ x_k,x_i^{\mp1},x_i,x_j \} \otimes e_k^*$,
we obtain
\begin{gather}
h_{ij} \otimes e_j^* \d h_{kj} \otimes e_j^* - h_{ki} \otimes e_i^*, \label{eq12} \\
h_{ij} \otimes e_i^* \d h_{ki} \otimes e_k^* - h_{kj} \otimes e_j^*, \label{eq14} \\
h_{ij} \otimes e_i^* \d h_{kj} \otimes e_k^* + h_{ki} \otimes e_i^*. \label{eq13}
\end{gather}

Now we show that all $h_{ij} \otimes e_j^*$ and $h_{ij} \otimes e_i^*$ are rewritten as a linear combination of $h_{1j} \otimes e_j^*$.
First, if $i, j \neq 1$, then $h_{ij} \otimes e_j^* \d h_{1j} \otimes e_j^* - h_{1i} \otimes e_i^*$ by (\ref{eq12}).
From (\ref{eq11}), we see $h_{ij} \otimes e_j^* \d h_{ji} \otimes e_i^*$, and hence $h_{i1} \otimes e_1^* \d h_{1i} \otimes e_i^*$
for any $i \neq 1$.

On the other hand, by (\ref{eq14}) and (\ref{eq13}),
\[ 2 h_{ij} \otimes e_i^* \d h_{ki} \otimes e_k^* - h_{kj} \otimes e_j^* + h_{kj} \otimes e_k^* + h_{ki} \otimes e_i^*. \]
Using (\ref{eq11}), (\ref{eq12}) and the results above, we obtain
\[\begin{split}
   2 h_{ij} \otimes e_i^* & \d - h_{ij} \otimes e_j^* - h_{kj} \otimes e_j^* + h_{ki} \otimes e_i^* \\
   & = -2 h_{ij} \otimes e_j^* \\
   & = \begin{cases}
        - 2 h_{1j} \otimes e_j^* + 2 h_{1i} \otimes e_i^*, \hspace{1em} & \mathrm{if} \hspace{1em} i,j \neq 1, \\
        - 2 h_{1i} \otimes e_i^*, & \mathrm{if} \hspace{1em} j=1.
       \end{cases}
  \end{split}\]

\vspace{1em}

From the argument above, we conclude that the generating set $\mathfrak{E}^*$ of
$\bar{R}^{\mathrm{ab}} {\otimes}_{\mathrm{Aut}^{+}F_n} H_L^*$ can be reduced to
\[\begin{split}
  \{ r_{i^{\pm1}j}(\cdot) \otimes e_p \,|\, & p \neq j \} \, \cup \, \{ r_{i^{\pm1}j}(1) \otimes e_j^* \,|\, i,j \neq 1 \} \\
  & \cup \, \{ r_{1^{\pm1}j}(2) \otimes e_j^* \,|\, j \neq 1,2 \} \cup \, \{ h_{1j} \otimes e_j^* \,|\, j \neq 1\}.
  \end{split}\]
The number of the generators above is just $2n(n^2-n)-n-1$.
Hence it is a basis of $\bar{R}^{\mathrm{ab}} {\otimes}_{\mathrm{Aut}^{+}F_n} H_L^*$
as a free $L$-module. This completes the proof of Proposition {\rmfamily \ref{P2}}. $\square$

\section{Acknowledgments}\label{6}

The author would like to thank Professor Nariya Kawazumi for
valuable advice and warm encouragement.
This research is supported by a JSPS Research Fellowships for Young Scientists.


\begin{thebibliography}{99}

 \bibitem{Bro} K. S. Brown; Cohomology of groups, Springer-Verlag, 1982.
 \bibitem{Ger} S. M. Gersten; A presentation for the special automorphism group of a free group,
               J. Pure and Applied Algebra 33 (1984), 269-279.
 \bibitem{Hat} A. Hatcher and K.\ Vogtmann; Rational homology of $\mathrm{Aut}(F_n)$, Math. Res.
               Lett. 5 (1998), 759-780.
 \bibitem{Nat} A. Hatcher and N. Wahl; Stabilization for the automorphisms of free groups with boundaries,
               Geometry and Topology, Vol. 9 (2005), 1295-1336.
 \bibitem{Hil} P. J. Hilton and U. Stammbach; A Course in Homological Algebra, Springer-Verlag, 1971.
 \bibitem{Kaw} N. Kawazumi; Cohomological aspects of Magnus expansions, preprint, The University of Tokyo. UTMS 2005-18 (2005).
 \bibitem{Kor} M. Korkmaz and A. Stipsicz; The second homology groups of mapping class groups of orientable surfaces,
               Math. Proc. Camb. Phil. Soc., 134 (2003), 479-489.
 \bibitem{Sat} T. Satoh; Twisted first homology group of the automorphism group of a free group,
               Journal of Pure and Applied Algebra 204 (2006), 334-348, to appear.
 \bibitem{Sa2} T. Satoh; Twisted second homology group of the automorphism group of a free group,
               preprint,  UTMS 2006-1 (2006),
               \texttt{http://kyokan.ms.u-tokyo.ac.jp/users/preprint/pdf/2006-1.pdf}.
\end{thebibliography}
\end{document}